%
%

\input amstex
\documentstyle{amsppt}
\NoBlackBoxes

\define\france{1}
\define\langue{\france}
\define\titre#1 #2 #3{{#2} {\kern-0.2em\bf (#3)} {#1}} 

\newcount\numbib\numbib=0\define\bib#1{{\global\advance\numbib by 1}\edef#1{\the\numbib}}

\bib\BUZt\bib\BUZb\bib\BUZf\bib\BUZk
\bib\FELa
\bib\FERa\bib\FERb\bib\FERc
\bib\GERa
\bib\GHAa
\bib\HOFb
\bib\KINa
\bib\NEWc
\bib\ORWa
\bib\OWb
\bib\PETa
\bib\ROTa
\bib\RUDa
\bib\STa

\define\thmLien{B}
\define\thmCaract{C}
\define\remPnonPistante{(1.9)}
\define\corCritere{(3.1)}
\define\defDistSyst{section~3.3}
\define\secNonPiste{5}

\def\picture #1 by #2 (#3){
  \vbox to #2{
    \hrule width #1 height 0pt depth 0pt
    \vfill
    \special{picture #3} 
    }
  }
\def\scaledpicture #1 by #2 (#3 scaled #4){{
  \dimen0=#1 \dimen1=#2
  \divide\dimen0 by 1000 \multiply\dimen0 by #4
  \divide\dimen1 by 1000 \multiply\dimen1 by #4
  \picture \dimen0 by \dimen1 (#3 scaled #4)}
  }

\define\eps{\epsilon}\define\vareps{\varepsilon}
\define\fb{\bar{f}}

\define\omegat{{\tilde{\omega}}}

\define\Z{\Bbb Z}\define\N{\Bbb N}

\define\diffsym{\operatorname{\Delta}}
\define\db{{\bar{d}}}

\define\eqdef{\overset{\operatorname{def}}\to{=}}

\define\card{\#}

\redefine\top{{\operatorname{top}}}

\define\new#1{{\bf #1}}
\define\stress#1{{\sl #1}}
\define\cte{\text{const}}

\def\today{
\if\langue\france
     \ifnum\day=1
     	 {1er\space}
     \else
       {\number\day\space}
     \fi
     \ifcase\month\or
        Janvier\or F\'evrier\or Mars\or Avril\or Mai\or Juin\or Juillet\or
        Ao\^ut\or Septembre\or Octobre\or Novembre\or D\'ecembre\fi\space
     \number\year
\else
     \ifcase\month\or
        January\or February\or March\or April\or May\or June\or July\or
        August\or September\or October\or November\or December\fi\space
\number\day,\space\number\year
\fi}

\define\section{*}
\newcount\numsec
\newcount\numsub
\newcount\numitm
\numsec=0\numsub=0
\define\newsec{
  \global\advance\numsec by 1
  \global\redefine\section{\the\numsec}\section. 
  \global\numsub=0\global\numitm=0
  }

\define\newsub{
  \global\advance\numsub by 1
  {\rm \section.\the\numsub.}
  }

\define\thethm{\section.\the\numitm}
\define\thm{
\if\langue\france
     \titre Th\'eor\`eme {\global\advance\numitm by 1} \thethm
\else
     \titre Theorem {\global\advance\numitm by 1} \thethm
\fi}

\define\thedef{\section.\the\numitm}
\define\defn{
\if\langue\france
     \titre D\'efinition {\global\advance\numitm by 1} \thedef
\else
     \titre Definition {\global\advance\numitm by 1} \thedef
\fi}

\define\theprop{\section.\the\numitm}
\define\prop{\titre Proposition {\global\advance\numitm by 1} \theprop}

\define\thelem{\section.\the\numitm}
\define\lem{
\if\langue\france
     \titre Lemme {\global\advance\numitm by 1} \thelem
\else
     \titre Lemma {\global\advance\numitm by 1} \thelem
\fi}

\define\thecor{\section.\the\numitm}
\define\cor{
\if\langue\france
     \titre Corollaire {\global\advance\numitm by 1} \thecor
\else
     \titre Corollary {\global\advance\numitm by 1} \thecor
\fi}

\define\theconj{\section.\the\numitm}
\define\conj{
\if\langue\france
     \titre Conjecture {\global\advance\numitm by 1} \theconj
\else
     \titre Conjecture {\global\advance\numitm by 1} \theconj
\fi}

\define\theclm{\section.\the\numitm}
\define\clm{
\if\langue\france
     \titre Affirmation {\global\advance\numitm by 1} \theclm
\else
     \titre Claim {\global\advance\numitm by 1} \theclm
\fi}

\define\therem{\section.\the\numitm}
\define\rem{
\if\langue\france
     \titre Remarque {\global\advance\numitm by 1} \therem
\else
     \titre Remark {\global\advance\numitm by 1} \therem
\fi}
\define\endrem{\qed\medbreak}

\define\theex{\section.\the\numitm}
\define\exemple{
\if\langue\france
     \titre Exemple {\global\advance\numitm by 1} \theex
\else
     \titre Example {\global\advance\numitm by 1} \theex
\fi}

\define\theequ{(\section.\the\numitm)}
\define\equ{{\global\advance\numitm by 1}\section.\the\numitm}

\define\testpredef#1#2{ 
    \edef\xxx{\the#2} 
    \ifx\xxx#1 \message{ok pour \noexpand#1}
       \else   \message{Mettez a jour \noexpand#1}
    \fi}
\define\testnumthm#1{\testpredef#1\numthm}
\define\testnumprop#1{\testpredef#1\numprop}

\topmatter
\title
    Sur une g\'en\'eralisation de la notion de syst\`eme dynamique de rang un
   d\'efinie par une propri\'et\'e de pistage${}^*$
\\
\\
\\
{\it (On~a~weak~version~of~the~rank~one~property defined~by~shadowing) }
\\
\\
J\'er\^ome BUZZI
\endtitle

\rightheadtext{Une g\'en\'eralisation du rang un} 



\abstract
We investigate a shadowing property which appears naturally in the study of piecewise monotonic maps of the interval. It
turns out to be a weak form of the rank one property, a well-known notion in abstract ergodic theory. We show that this
new property is implied by finite or even local rank, but that it is logically independent of loose Bernoulliness. We
give (counter)examples, including L.B.\ systems with arbitrarily high-order polynomial complexity. The shadowing
property defines a small subset of all zero-entropy systems, in the sense that it defines a closed set with empty
interior with respect to the $\db$-metric, induced by the Hamming distance.

We also make some remarks on the link between the shadowed system and the sequence assumed by
the shadowing property.

\vskip1cm

\noindent
{\smc R\'esum\'e.}
Nous \'etudions une propri\'et\'e de pistage qui appara\^{\i}t naturellement dans l'\'etude des applications monotones
par morceaux sur l'intervalle. Cette propri\'et\'e s'av\`ere \^etre un affaiblissement de la notion de rang un, notion
bien connue en th\'eorie ergodique abstraite. Nous montrons que cette nouvelle propri\'et\'e est impliqu\'ee par le rang
fini ou m\^eme local, mais qu'elle est logiquement ind\'ependante de la l\^ache Bernoullicit\'e. Nous donnons des
(contre)exemples, en particulier un syst\`eme L.B.\ ayant une complexit\'e d'ordre polyn\^omial arbitrairement
\'elev\'e. La propri\'et\'e de pistage d\'efinit une petite partie de tous les syst\`emes d'entropie nulle en ce qu'il
s'agit d'un ferm\'e d'int\'erieur vide au sens de la m\'etrique $\db$, induite par la distance de Hamming.

Nous faisons \'egalement quelques remarques sur le lien entre le syst\`eme dynamique pist\'e et la suite pistante. 

\endabstract





\endtopmatter

\document

\noindent
${}^*$ Travail en partie effectu\'e au Laboratoire de Topologie de Dijon, UMR 5584.
\vskip1cm

\newpage

\head \newsec Introduction \endhead

Cet article est consacr\'e \`a l'analyse abstraite de la propri\'et\'e de pistage d'un syst\`eme dynamique
probabiliste par une petite partie de l'espace sous-jacent (cf.\ les d\'efinitions ci-dessous). Cette
propri\'et\'e appara\^{\i}t naturellement dans l'\'etude de certains sys\-t\`e\-mes dynamiques
g\'eo\-m\'e\-tri\-ques avec singularit\'es: elle distingue les mesures invariantes ``d\'eg\'en\'er\'ees'',
c'est-\`a-dire dont les points g\'en\'eriques se rapprochent trop vite de ces singularit\'es. Il est
donc important de comprendre quelles contraintes cette propri\'et\'e impose \`a ces mesures
d\'eg\'en\'er\'ees, voire \`a leur existence.

On se concentre ici sur le cas des syst\`emes $1$-pist\'es c'est-\`a-dire pist\'es par une partie
r\'eduite \`a un point. Ce cas correspond g\'eom\'etriquement aux mesures d\'eg\'en\'er\'ees
des applications monotones par morceaux sur l'intervalle. On montre que les syst\`emes $1$-pist\'es
forment une nouvelle classe de syst\`emes dynamiques d'entropie nulle. Cette classe g\'en\'eralise
les syst\`emes de rang fini ou local. On donne une caract\'erisation des syst\`emes $1$-pist\'es calqu\'ee
sur celle des syst\`emes de rang un ou fini. 

Une partie importante de nos efforts est consacr\'ee \`a la construction d'exemples qui montrent en particulier que le
$1$-pistage est ind\'ependant de la propri\'et\'e ``l\^a\-che\-ment Bernoulli''. Enfin on envisage les liens entre
l'orbite du point pistant et le syst\`eme pist\'e.

\subheading{\newsub D\'efinition du pistage}
On part de la propri\'et\'e symbolique suivante. Soit $\Cal A$ un ensemble fini et $\Sigma$ une partie
quelconque de $\Cal A^\N$. Une suite $A\in\Cal A^\Z$ est dite \new{$\Sigma$-pist\'ee} s'il
existe des entiers $n_i,m_i$ tendant vers l'infini et une suite d'\'el\'ements $s^{(i)}\in\Sigma$ tels
que:
 $$
           A_{-n_i}\dots A_{+m_i} = s^{(i)}_0\dots s^{(i)}_{n_i+m_i}.
 $$

Soit $(X,\Cal B,T,\mu)$ un syst\`eme dynamique probabiliste, c'est-\`a-dire un automorphisme d'un
espace de Lebesgue. Soit $P$ une partition de $X$ (suppos\'ee finie et mesurable comme toutes les
partitions de cet article)  et $\Sigma$ une partie quelconque de $P^\N$. Pour $F\subset\Z$, on
d\'efinit le $P,F$-nom de $x$, $P^F(x)$, comme l'application $A$ de $F$ dans $P$ d\'efinie par
$f^k(x)\in A_k$ pour $k\in F$.

\proclaim{\defn}
On dira que $(T,\mu)$ est \new{$(P,\Sigma)$-pist\'e} si $P$ est une partition g\'en\'eratrice et si, pour $\mu$-presque
tout $x\in X$, le $P$-nom de $x$, $P^\Z(x)$, est $\Sigma$-pist\'e.

On dira que $(T,\mu)$ est \new{$1$-pist\'e} (ou simplement, \new{pist\'e}) s'il est $(P,\Sigma)$-pist\'e pour une
certaine partition $P$ et une partie $\Sigma$ r\'eduite \`a un \'el\'ement $\omega\in P^\N$. On dit que $P$ et
$\omega$ sont, respectivement, une partition et une suite \new{pistantes}.
 \endproclaim

\remark{\rem}
Les propri\'et\'es: pist\'e par $\Sigma$ fini et pist\'e par $\Sigma$ r\'eduit \`a un \'el\'ement sont clairement
\'equivalentes pour les syst\`emes ergodiques. 
 \endrem\endremark

L'int\'er\^et de cette notion est tout d'abord qu'elle {\bf appara\^{\i}t naturellement} dans certains sys\`emes
g\'eom\'etriques. Donnons un exemple particuli\`erement simple:

\remark{\exemple} \edef\exEch{(\theex)}
Un \'echange d'intervalles $([0,1],T)$ muni d'une mesure invariante et ergodique est $1$-pist\'e. 
\endrem\endremark

C'est un cas particulier du fait que les mesures ``d\'eg\'en\'er\'ees au sens de Hofbauer'' d'une application monotone
par morceaux de l'intervalle sont exactement les mesures pist\'ees par les itin\'eraires \`a droite et \`a gauche des
points critiques (i.e., les  ``kneading invariants''). Par mesure d\'eg\'en\'er\'ee, on entend ici une mesure port\'ee
par la ``partie non-markovienne'' de l'extension naturelle de la dynamique symbolique (cf.\ \cite{\HOFb}, 
\cite{\NEWc}). Cette situation s'\'etend aux g\'en\'eralisations de la construction de F.~Hofbauer (\cite{\BUZt,
chap.~6}, \cite{\BUZb, pp. 145--149}, \cite{\BUZf}).

Le pistage joue \'egalement un r\^ole dans la construction des vari\'et\'es stables ou instables des
syst\`emes avec singularit\'es. Expliciter ce r\^ole est parfois fructueux \cite{\BUZk}.

\demo{Preuve directe}
On traite le cas topologiquement minimal, le cas g\'en\'eral s'en d\'e\-dui\-sant facilement.
On va montrer que le $P$-nom de tout point est $\Sigma$-pist\'e avec $P$ la partition naturelle en intervalles et
$\Sigma$ les itin\'eraires \`a droite et \`a gauche des discontinuit\'es. Notons $D$ l'ensemble de ces discontinuit\'es.
Soit $x\in[0,1]$. On laisse au lecteur le cas exceptionnel o\`u $x$ est sur l'orbite d'une discontinuit\'e.

Par minimalit\'e,  $\liminf_{n\to\infty} d(T^{-n}x,D)=0$.  Il existe donc $n_1<n_2<\dots$ tels que
$d(T^{-n_i}x,D)<d(T^{-n}x,D)$ pour $0\leq n<n_i$. Soit $d_i\in D$ r\'ealisant la distance pr\'ec\'edente. Posons $\pm=+$
si $d_i<x$, $\pm=-$ sinon. $T$ \'etant une isom\'etrie sur chacun de ses intervalles, l'in\'egalit\'e stricte
pr\'ec\'edente implique que pour $0\leq k\leq n_i$, $T^{-n_i+k}x$ et $T^kd_i\pm$ ne sont pas s\'epar\'es par une
discontinuit\'e: ils sont donc dans le m\^eme \'el\'ement de $P$. C'est pr\'ecis\'ement le pistage annonc\'e.
 \qed \enddemo

\remark{\rem}
On verra \remPnonPistante\ qu'une translation irrationnelle sur $\Bbb T^d$ muni de la mesure de Lebesgue est \'egalement
$1$-pist\'ee mais que, pour $d\geq2$, l'\'equivalent de la partition naturelle n'est pas, dans certains cas, admissible
pour le pistage. 
\endrem\endremark

\medbreak

Une autre motivation \`a cette d\'efinition c'est qu'il s'agit d'une {\bf g\'en\'eralisation naturelle} de la
notion de rang un. D'apr\`es le th\'eor\`eme~\thmCaract\ ci-dessous, on a en effet la:

\proclaim{Caract\'erisation}
$(T,\mu)$ est $1$-pist\'e si et seulement si pour tout $\eps>0$ et toute partition $Q$ il existe une tour irr\'eguli\`ere
$\eps$-raffinant $Q$, i.e., un ensemble mesurable $B$ tel que si $h:B\to\N^*$ est le temps de retour dans $B$, alors,
pour tout $A\in Q$, il existe une union $A'$ d'ensembles de la forme: $T^k\{x\in B:h(x)>k\}$, avec la propri\'et\'e:
$\sum_{A\in P} \mu(A\diffsym A')<\eps$.
\endproclaim

Remarquons que cette caract\'erisation ne fait intervenir ni partition ni suite pistantes.

\remark{\rem}
On peut donner une d\'efinition {\sl m\'etrique} du pistage. Dans le cas du $1$-pistage,
l'\'equivalence des deux notions se d\'eduit du th\'eor\`eme~\thmCaract\ ci-dessous.
\endrem
\endremark

\subheading{\newsub R\'esultats}

Rappelons tout d'abord le fait suivant:

\proclaim{\thm\ (\cite{\BUZt, (6.4)} ou \cite{\BUZb, Theorem~6.1})}
Si $(T,\mu)$ est $(P,\Sigma)$-pist\'e alors:
 $$
    h(T,\mu) \leq h_\top(\Sigma,\sigma) \eqdef \limsup_{n\to\infty} \frac1n\log\card\{s_0\dots s_{n-1}:s\in\Sigma\}.
 $$
\endproclaim

Dans le cas $1$-pist\'e, l'entropie topologique $h_\top(\Sigma,\sigma)$ est bien \'evidemment {\bf nulle} (et ceci m\^eme
si l'orbite pistante est dense): la classe des syst\`emes $1$-pist\'es est donc une partie de l'ensemble des syst\`emes
d'entropie nulle.

Evidemment, dire qu'une partie est r\'eduite \`a un point est bien plus fort que de dire que son entropie est
nulle. On peut donc esp\'erer que le $1$-pistage soit une notion bien plus forte que l'entropie nulle: c'est la
motivation de ce travail.

\subheading{Position de la classe des syst\`emes $1$-pist\'es}
Notre r\'esultat principal situe les syst\`emes $1$-pist\'es par rapport \`a des classes de syst\`emes d'en\-tro\-pie
nulle bien connues en th\'eorie ergodique abstraite (cf.\ \cite{\FERa,\FERc,\KINa,\ORWa}):

\proclaim{\thm\ A} 
L'ensemble des syst\`emes $1$-pist\'es est une partie de l'ensemble des syst\`emes d'entropie nulle. Au sens de la
distance de Hamming (\defDistSyst), c'est une partie ferm\'ee d'int\'erieur vide.

On a les inclusions suivantes (chaque propri\'et\'e d\'esignant l'ensemble des syst\`emes
dynamiques probabilistes ergodiques la v\'erifiant):
 $$\gather
      \text{rang local}
      \subsetneq \text{$1$-pistage} \subsetneq \text{entropie nulle}\\
     \text{$1$-pistage} \not\subset \text{l\^achement Bernoulli}
    \text{ et } \text{l\^achement Bernoulli} \not\subset \text{$1$-pistage}.
 \endgather$$
 \endproclaim
(on sait que: rang un $\subsetneq$ rang fini $\subsetneq$ rang local).

On en d\'eduit que les rotations sur les groupes compacts, les \'echanges d'in\-ter\-val\-les ou encore les substitutions
sont $1$-pist\'es (ce sont des syst\`emes de rang un ou fini \cite{\FERc}).

\medbreak

{\smc Probl\`eme:} Trouver un exemple naturel de syst\`eme d'entropie nulle qui ne soit pas $1$-pist\'e.
\medbreak

{\smc Probl\`eme:} Trouver une classe naturelle de syst\`emes dynamiques g\'eom\'etriques qui soient $1$-pist\'es mais
dont certains ne soient pas de rang local.

(que peut-on dire des \'echanges isom\'etriques de polygones \cite{\GHAa} ?)
\medbreak

{\smc Probl\`eme:} Trouver des propri\'et\'es ergodiques impliqu\'ees par le $1$-pistage (en dehors de l'entropie
nulle) et notamment: 

 ---le $1$-pistage implique-t-il, comme le rang local, que la multiplicit\'e spectrale est finie?

 ---implique-t-il une certaine vitesse de r\'ecurrence en analogie avec les r\'esultats de
D.~Ornstein et B.~Weiss \cite{\OWb} qui disent que, pour $\mu$-presque tout $x$:
 $$
        \lim_{n\to\infty} \frac1n\log \min\{j\geq n:P^{[j,j+n[}(x)=P^{[0,n[}(x) \} = h_\mu(T,P).
 $$

\comment
\subheading{Stabilit\'e du $1$-pistage \rm (par rapport aux proc\'ed\'es de constructions usuels)}
\proclaim{\prop} 
Un {\bf it\'er\'e}, un {\bf facteur} ou une {\bf limite inductive} de syst\`emes $1$-pist\'es sont encore $1$-pist\'es. 

Au contraire une {\bf suspension} de syst\`eme $1$-pist\'e n'est pas n\'ecessairement $1$-pist\'ee.
\endproclaim

 \proclaim{\conj}
Les op\'erations de produit cart\'esien, d'extension finie et d'in\-duc\-tion sur un ensemble
mesurable peuvent \'egalement d\'etruire le $1$-pistage. Il est m\^eme vraisemblable que, un
syst\`eme $1$-pist\'e non-p\'eriodique et ergodique \'etant fix\'e, l'in\-duc\-tion par presque tout
mesurable d\'etruit le $1$-pistage. 
 \endproclaim
\endcomment

L'existence de syst\`emes d'entropie nulle non $1$-pist\'es, ayant de plus la propri\'et\'e L.B.\ (dont la d\'efinition
est rappel\'ee au d\'ebut de la section~4) se d\'eduit du r\'esultat suivant:

\proclaim{Th\'eor\`eme (5.1)}
Pour tout $\Gamma<\infty$, il existe $(X,\Cal B,T,\mu)$ un automorphisme ergodique d'un espace de
Lebesgue qui est d'entropie nulle, l\^achement Bernoulli et qui admet une partition g\'en\'eratrice
$P$ telle que, pour $\eps>0$ assez petit, pour tout $n$ assez grand, la mesure de toute boule-$\db$ (d\'efinie au
d\'ebut de la section~5) de rayon $\eps$ correspondant \`a un mot de longueur $n$ est major\'ee par $n^{-\Gamma}$.
\endproclaim

\subheading{Lien entre suite pistante et syst\`eme pist\'e}
Presque toute orbite \'etant form\'ee de copies exactes de d\'ebuts arbitrairement longs de la suite pistante on
pourrait penser que ce lien est tr\`es fort. Nous pr\'esentons quelques observations qui
montrent qu'en g\'en\'eral il n'en est rien. 

Rappelons qu'un point $\omega$ de l'espace topologique $P^\N$ est \new{quasi-g\'en\'erique} pour une mesure $\mu$, si 
la suite $\frac1n\sum_{k=0}^{n-1} \delta_{\sigma^k\omega}$ admet $\mu$ comme valeur d'adh\'erence pour la topologie
vague. S'il y a convergence vers $\mu$, alors $\omega$ est dit \new{g\'en\'erique}.

\proclaim{\thm\ \thmLien}
Soit $(X,\Cal B,T,\mu)$ un syst\`eme dynamique probabiliste ergodique.

\medbreak

(1)~Si $\omega\in P^\N$ piste $(T,\mu)$ par rapport \`a une certaine partition $P$ alors $\omega$ est
quasi-g\'en\'erique, mais non-n\'ecessairement g\'en\'erique, pour $P^\Z\mu$.

\medbreak

(2)~Il y a abondance de suites pistantes:

Au sens de la {\bf topologie,} presque toute suite piste au moins un syst\`eme dynamique. Plus
pr\'ecis\'ement, si $P$ est une partition de $X$, alors l'en\-sem\-ble des $\omega\in P^\N$ tels que $(T,\mu)$ soit
$(P,\omega)$-pist\'e pour au moins une mesure $\mu$ est un $G_\delta$-dense de $P^\N$.

Au sens de la {\bf mesure,} si $T$ est de rang local et si $P$ est une partition par rapport \`a laquelle $T$
soit pist\'ee par une certaine suite alors le $P$-nom, $P^\Z(x)$, de $\mu$-presque tout $x\in X$, $P$-piste $T$.

\medbreak

(3)~Une m\^eme suite peut pister des syst\`emes compl\`etement diff\'erents: 

Etant donn\'es deux syst\`emes $1$-pist\'es, ergodiques et ap\'eriodiques (i.e., dont les points p\'e\-rio\-di\-ques
forment un ensemble de mesure nulle), on peut trouver une partition g\'en\'eratrice \`a deux \'el\'ements pour chacun des
syst\`emes, telles que chaque syst\`eme est pist\'e, par rapport \`a la partition choisie, par une m\^eme suite $\omega$
(si on identifie les deux partitions).
 \endproclaim

La derni\`ere partie du th\'eor\`eme sugg\`ere la:

{\smc Question:} Existe-t-il une suite {\bf universellement pistante} $\Omega\in\{0,1\}^\N$, c'est-\`a-dire telle que
tout syst\'eme ergodique et ap\'eriodique admette une partition g\'e\-n\'e\-ra\-tri\-ce $\{A_0,A_1\}$ par rapport \`a
laquelle ce syst\`eme soit pist\'e par la suite $A_{\Omega_0}A_{\Omega_1}\dots$?
\medbreak

On peut interpr\'eter cette derni\`ere partie du th\'eor\`eme comme indiquant, qu'{\bf en g\'en\'eral},
la connaissance d'une suite pistante {\bf n'apprend rien} sur le syst\`eme $1$-pist\'e. Ceci sugg\`ere donc de regarder
des suites pistantes particuli\`eres. On pose donc le:
\medbreak

{\smc Probl\`eme:} Formuler et prouver des th\'eor\`emes du type: une propri\'et\'e ``combinatoire'' de la suite
pistante implique l'existence ou l'unicit\'e ou une propri\'et\'e ergodique pour le syst\`eme pist\'e. 

On aimerait tout particuli\`erement:
 \roster
  \item des conditions sur la suite qui puissent \^etre v\'erifi\'ees dans les cas
g\'eo\-m\'e\-tri\-ques, notamment pour les itin\'eraires des points critiques des applications
monotones par morceaux.
  \item identifier ``combinatoirement'' les sous-classes rang un, rang fini, rang local parmi les syst\`emes
pist\'es.
 \endroster

\subheading{\newsub Choix de la partition}
Toute partition n'est pas pistante. Donnons un exemple simple.

\remark{\exemple}
Une rotation irrationnelle sur le cercle est un \'echange d'intervalles minimal: d'apr\`es l'exemple~\exEch, ce syst\`eme
est pist\'e par rapport \`a sa partition naturelle. C'est faux en dimension sup\'erieure: une translation
irrationnelle sur $\Bbb T^2$, qu'on peut voir comme un produit de deux rotations irrationnelles admet une partition
naturelle $P$ telle que les atomes des it\'er\'es sont de mesure (de Lebesgue) en $1/n^2$ si les angles sont \`a
quotients partiels born\'es. Or $\sum 1/n^2<\infty$: on voit alors facilement (cf.\ par exemple \corCritere) que le
syst\`eme ne peut \^etre pist\'e par rapport \`a sa partition naturelle!
\endrem\endremark

\subheading{\newsub Plan de l'article}
Tout d'abord, la section~2 donne des d\'efinitions alternatives du $1$-pistage calqu\'ees sur les caract\'erisations
classiques du rang $1$. La section~3 en tire quelques cons\'equences qui prouvent, d'une part, les
deux premiers points du th\'eor\`eme~A, d'autre part, le th\'eor\`eme~\thmLien.

\medbreak

Le reste de l'article donne deux exemples achevant la preuve du th\'eor\`eme~A. La section~4 d\'ecrit un syst\`eme
$1$-pist\'e qui n'est pas L.B. Sa construction est inspir\'ee d'un exemple de syst\`eme d'entropie z\'ero non-L.B.\ %
d\^u \`a J.~Feldman. Ici toutefois il y a un type de bloc pr\'edominant (affect\'e d'une probabilit\'e $1/n$ au rang
$n$) qui force le $1$-pistage mais rend plus d\'elicate la preuve du caract\`ere non-L.B. La section~5 d\'ecrit
inversement un syst\`eme L.B., d'entropie nulle qui n'est pas $1$-pist\'e.  Pour ce faire, on construit des syst\`emes
L.B.\ dont la complexit\'e mesur\'ee est d'ordre polynomial arbitrairement \'elev\'e, en codant le type de chaque bloc
par de petits d\'ecalages (importants pour la m\'etrique $\db$, n\'egligeables pour $\fb$, d'o\`u le r\'esultat).

\medbreak

\subheading{\newsub Remerciements}
Je dois \`a Philippe {\smc Thieullen} l'id\'ee d'un lien possible entre $1$-pistage et rang un.

Je remercie Jean-Paul {\smc Thouvenot} pour de nombreuses discussions auxquelles ce travail doit beaucoup.

\head \newsec  Le pistage comme rang $1$ irr\'egulier \endhead

Pour pr\'eciser la relation entre le $1$-pistage et la notion de rang un et ses
g\'e\-n\'e\-ra\-li\-sa\-tions nous montrons qu'on peut caract\'eriser le $1$-pistage d'une fa\c{c}on
similaire au rang un. Ce faisant nous obtiendront une certaine souplesse vis-\`a-vis de la suite pistante, souplesse
qui nous servira par la suite. 

Nous avons besoin de quelques d\'efinitions.

 \proclaim{\defn}
Une \new{tour irr\'eguli\`ere} $\tau$ est d\'efinie par la donn\'ee d'une fonction mesurable et
born\'ee, sa \new{hauteur}, $h_\tau:X\to\N$ telle que les ensembles suivants, appel\'es
\new{niveaux}, soient deux-\`a-deux disjoints:
 $$
        N_\tau(k) = T^k\{x\in X: h_\tau(x)>k\} \qquad k=0,1,\dots,\max h_\tau-1.
 $$
L'union des niveaux d'une tour est appel\'e son \new{support} et not\'e simplement $\tau$.
Le \new{to\^{\i}t} est l'ensemble $\bigcup_{k\geq1} T^{k-1}(h_\tau^{-1}(k))$, i.e., l'ensemble des $x\in\tau$ qui
``sortent'' de la tour. 
La \new{partition associ\'ee}, not\'ee $P_\tau$, est la partition
engendr\'ee par les niveaux de $\tau$, i.e., les $N_\tau(k)$ pour $0\leq k<\max h_\tau$.

Etant donn\'ee une partition $P$, une tour est dite \new{$P$-pure} si chacun de ses niveaux est inclus
dans un \'el\'ement de $P$.
 \endproclaim

 \remark{\rem}
1.~La condition de disjonction ci-dessus est \'equivalente au fait que la fonction hauteur est major\'ee par le temps de
retour \`a la base $\{x: h_\tau(x)>0\}$.

2.~Si la fonction hauteur ne prend qu'une seule valeur non-nulle $H$, on obtient simplement une {\bf tour de Rokhlin}
de hauteur $H$ et de base $h_\tau^{-1}(H)$. 
 \endrem\endremark

\proclaim{\defn}
Une tour irr\'eguli\`ere $\tau'$ est dite \new{embo\^{\i}t\'ee} dans une tour de Rokhlin de base $B$ et de hauteur $H$
s'il existe des entiers $l_1<\dots<l_r$ tels que:
 $$
     T^k B = \bigcup_{i=1}^r N_{\tau'}(l_i+k) \qquad \forall 0\leq k<H.
 \tag\equ$$ \edef\equRokhEmb{\theequ}

$\tau$ et $\tau'$ \'etant deux tours irr\'eguli\`eres, on dit que $\tau'$ est \new{embo\^{\i}t\'ee}
dans $\tau$ si $\tau'$ s'embo\^{\i}te dans chaque \new{sous-tour de Rokhlin} de $\tau$. Les
\new{sous-tours de Rokhlin} de $\tau$ sont les tours de base $h_\tau^{-1}(k)$ et de hauteur $k$, $k$ parcourant les
hauteurs de $\tau$.
 \endproclaim

\equRokhEmb\ est \'equivalent \`a: $B=\bigcup_{i=1}^r N_{\tau'}(l_i)$ et $]l_i,l_i+H[$ ne contient pas de
\new{hauteurs} de $\tau'$, i.e., d'entiers $k>0$ tels que $h_{\tau'}^{-1}(k)\ne\emptyset$.

On peut \'egalement formuler l'embo\^\i tement de $\tau'$ dans $\tau$ en disant que $P_\tau'$ est plus fine que
$P_\tau$ et que $\tau'\setminus\text{to\^\i t}(\tau')$ contient $\tau\setminus\text{to\^\i t}(\tau)$; ou encore que
l'ordre partiel naturellement associ\'e \`a $\tau'$ prolonge celui correspondant \`a $\tau$.
\medbreak

On note enfin $\db(a,b)\eqdef\frac1n\card\{k=0,\dots,n-1:a_k\ne b_k\}$ la distance de Hamming entre deux
suites finies de m\^eme longueur $n=|a|=|b|$.

On peut maintenant \'enoncer le:

\proclaim{\thm\ \thmCaract}
Le $1$-pistage d'un syst\`eme $(X,\Cal B,T,\mu)$ ergodique est caract\'eris\'e par n'importe
laquelle des propri\'et\'es suivantes:

(C1) pour toute partition $P$, il existe une suite $\omega\in P^\N$ telle que, pour
$\mu$-presque tout $x\in X$, il existe des entiers positifs $n_i$ et $m_i$ avec $n_i+m_i$ tendant vers l'infini tels
que:
 $$
        \lim_{i\to\infty} \db(P^{[-n_i,m_i]}(x),\omega_0\dots\omega_{n_i+m_i}) = 0.
 $$
 \smallbreak

(C2) pour toute partition $P$ et tout $\eps>0$, il existe une suite $\omega\in P^N$ ($N\in\N$),
telle que, pour $\mu$-presque tout $x\in X$, pour tout $n$ assez grand, la suite finie $P^{[0,n[}(x)$
peut s'\'ecrire comme la concat\'enation:
 $$
      \alpha^{(1)}\omega^{(1)}\alpha^{(2)}\omega^{(2)}\dots\alpha^{(r+1)}
 $$
o\`u les $\alpha^{(i)}$ sont des suites quelconques dont la somme des longueurs est au plus $\eps n$
et chaque $\omega^{(i)}$ v\'erifie $l=|\omega^{(i)}|\in[\eps^{-1},N]$ et $\db(\omega^{(i)},
\omega_0\dots\omega_{l-1})<\eps$. On dit qu'on a un \new{$(1-\eps)$-recouvrement} de $P^{[0,n[}(x)$ par
des \new{$\db$-$\eps$-copies} de d\'ebuts de $\omega$.
 \smallbreak

(C1') et (C2'): la m\^eme chose que (C1) et (C2), mais pour une partition g\'en\'eratrice
particuli\`ere.
 \smallbreak

(C3) pour toute partition $Q$ et tout $\eps>0$, il existe une tour irr\'eguli\`ere $\tau$ dont la
partition associ\'ee $P_\tau$ \new{$\eps$-raffine} $Q$, c'est-\`a-dire qu'\`a chaque \'el\'ement $A\in
Q$, on peut associer $A'$ une union d'\'el\'ements de $P_\tau$ de sorte que: $\sum_{A\in Q} \mu(A
\diffsym A')<\eps$.
 \smallbreak

(C4) il existe une suite de tours irr\'eguli\`eres toutes pures par rapport \`a une partition g\'en\'eratrice fix\'ee,
deux-\`a-deux em\-bo\^{\i}\-t\'ees et dont les partitions associ\'ees engendrent ensemble la tribu des mesurables modulo
$\mu$.
 \endproclaim

\remark{\rem}
Si on substitue, dans les caract\'erisations du th\'eor\`eme ci-dessus, des tours de
Rokhlin aux tours irr\'eguli\`eres, des mots de longueur fix\'ee $l=|\omega|$ aux segments initiaux
de longueur variable, {\bf on retrouve les caract\'erisations classiques des syst\`emes de rang un
\cite{\KINa}}. On voit d\`es \`a pr\'esent que le rang un implique le $1$-pistage (cf.\ section 3.1).
\endrem\endremark

\demo{Preuve du th\'eor\`eme~\thmCaract}
Il s'agit de prouver que les caract\'erisations propos\'ees par le th\'eor\`eme~\thmCaract\ sont
effectivement \'equivalentes au $1$-pistage. On va proc\'eder dans l'ordre suivant:
 \block
    $1$-pistage $\implies$ (C1) $\implies$ (C2) $\implies$ (C3) $\implies$ (C4) $\implies$ $1$-pistage.
 \endblock
Les \'equivalences de (C1)$\iff$(C1') et (C2)$\iff$(C2') seront envisag\'ees
s\'e\-pa\-r\'e\-ment.

Les d\'emonstrations s'inspirent, pour une large part, des d\'emonstrations ``classiques'' faites
dans le cas du rang un (voir notamment \cite{\KINa}).

\smallbreak

Remarquons que si $T$ n'est pas ap\'eriodique, alors, par ergodicit\'e, $X$ est r\'eduit \`a une
orbite p\'eriodique et le th\'eor\`eme est trivial. On suppose donc l'ap\'eriodicit\'e.

\smallbreak

\remark{\rem} \edef\remPurete{\therem}
Dans la cha\^{\i}ne d'implications $1$-pistage$\implies$(C1)$\implies$\dots$\implies$(C4), si $P$ est une
partition par rapport \`a laquelle le syst\`eme est pist\'e, alors:
 \roster
  \item dans (C1) et (C2), il n'y a pas d'erreur $\db$.
  \item dans (C3) et (C4), les tours obtenues peuvent \^etre suppos\'ees $P$-pures.
 \endroster
Cette remarque servira \`a la preuve de l'abondance en mesure des suites pistantes (section~3.2).
\endremark
\enddemo

\subheading{$1$-pistage $\implies$ (C1)}
Supposons donc $(X,\Cal B,T,\mu)$ $(Q,\omega)$-pist\'e avec $Q$ une partition g\'en\'eratrice.
Fixons $P$ une partition quelconque. $Q$ \'etant g\'en\'eratrice, il existe, pour chaque $\eps>0$,
un entier $N=N(\eps)$ tel que $P$ est $\eps$-raffin\'e par $\bigvee_{k=0}^{N-1} T^{-k}Q$ (un g\'en\'erateur est
automatiquement un g\'en\'erateur unilat\'eral, $T$ \'etant d'entropie nulle). Autrement dit, on a une application
$q_\eps:Q^N\to P$ telle que $p_\eps(Q^{[0,N[}(x))=P(x)$ pour tout $x\in X\setminus X_\eps$ avec $\mu(X_\eps)<\eps$. 

Posons $\eps_k=2^{-k}$, $N_k=N(\eps_k)$ et $X_k=X_{\eps_k}$. On d\'efinit $\omegat\in P^\N$ par:
 $$
       \omegat_n \eqdef p_{\eps_k}(\omega_n\dots\omega_{n+N_k-1}) \qquad \forall n=0,1,2,\dots
 $$
o\`u $k=k(n)\eqdef\max\{l:N_l\leq\log (n+1)+N_0\}$ (de sorte que $k(n)\nearrow\infty$ mais
$N_{k(n)}/n\to0$). On a: $\mu(\bigcup_{l\geq k} X_l)<2^{-k+1}$.

Soit $\eps>0$. Fixons $K<\infty$ tel que $2^{-K}<\eps/10$. D'apr\`es le th\'eor\`eme ergodique
appliqu\'e \`a $T^{-1}$, pour $\mu$-presque tout $x\in X$, si $n$ est
suffisamment grand:
 $$
        \frac{1}{n}\card\{p\in[-n,0]:T^p(x)\in \bigcup_{k\geq K} X_k\} < 2^{-K+1}
            < \frac{2}{10} \eps.
 \tag\equ$$\edef\equPetiteFrq{\theequ}
Fixons $n$ comme dans la d\'efinition du pistage: $Q^{[-n,0]}(x)=\omega_0\dots \omega_{n}$.
Posons $D=\frac{\eps}{10}n$. On peut supposer $n$ assez grand pour avoir l'in\'egalit\'e \equPetiteFrq\ ainsi
que $k(D)\geq K$ et $(\log(n+1)+N_0)/n<\eps/10$. En particulier, $N_{k(n)}\leq D$.

Soit $p\in[-n+D,-D]$ tel que $T^p(x)\notin \bigcup_{k\geq K} X_k$. Ces $p$ repr\'esentent une fraction de $[-n,0]$ au
moins \'egale \`a $1-4\eps/10$.

Comme $k=k(p+n)\leq k(n)$, $N_k\leq N_{k(n)}\leq D$ donc $p+N_k\leq 0$. En particulier, $\omega_{p+n}\dots
\omega_{p+n+N_k-1} =Q^{N_k}(T^px)$. Comme $k=k(p+n)\geq k(D)\geq K$ et $T^px\notin \bigcup_{l\geq
K} X_l$, 
 $$
     \omegat_{p+n} \eqdef p_{\eps_k}(\omega_{p+n}\dots\omega_{p+n+N_k-1})
                   = p_{\eps_k}(Q^{N_k}(T^px))=P(T^px). 
 $$
On en d\'eduit:
 $$
        \db(P^{[-n,0]}(x),\omegat_0\dots\omegat_n) < \frac{4}{10}\eps \leq \eps.
 $$
La propri\'et\'e (C1) est maintenant imm\'ediate.

\subheading{(C1)$\iff$(C1')}
L'implication (C1)$\implies$(C1') est triviale. L'implication r\'e\-ci\-pro\-que se d\'eduit de la preuve
ci-dessus.

\subheading{(C1) $\implies$ (C2)}
Il suffit d'appliquer le r\'esultat technique suivant, d\^u \`a D.~Rudol\-ph:

\proclaim{\thm\ (D.~Rudolph \cite{\RUDa, Theorem~3.6})}
(``Back\-ward Vi\-ta\-li Lem\-ma'')

Soit $(X,\Cal B,T,\mu)$ un automorphisme ap\'eriodique d'un espace de Lebesgue.

Supposons que pour presque tout $x$, on a d\'efini des entiers positifs $i_k$ et $j_k$ dont la somme
tend vers l'infini. Alors, pour tout $\eps>0$, il existe $Y\subset X$ et $k:Y\to\N$ tels que, si on pose $I(y)=\{T^ny :
-i_k\leq n\leq j_k\}$ alors:
 \roster
  \item Si $y,y'\in Y$ et $y\ne y'$ alors $I(y)\cap I(y')=\emptyset$.
  \item $\bigcup_{y\in Y} I(y)$ est de mesure au moins $1-\eps$.
 \endroster
\endproclaim

On applique le th\'eor\`eme avec $i_k(x)=n_k(x)$ et $j_k(x)=m_k(x)$ fournis par (C1). Quitte \`a r\'eduire un peu $Y$, on
peut supposer que $i_k+j_k$ est born\'ee. On applique ensuite le th\'eor\`eme ergodique pour obtenir (C2).

\subheading{(C2)$\iff$(C2')}
L'implication (C2)$\implies$(C2') est triviale.  L'implication r\'e\-ci\-pro\-que se d\'eduit de la
preuve de $1$-pistage$\implies$(C1) donn\'ee ci-dessus.

\subheading{(C2) $\implies$ (C3)}
Il suffit de reproduire la preuve correspondante dans le cas du rang un \cite{\KINa}. Rappelons
l'id\'ee. Soit $\eps>0$ et $Q$ une partition. 

$T$ \'etant ap\'eriodique, il existe une tour de Rokhlin de mesure presque $1$ et de grande hauteur
$H$. On subdivise sa base $B$ suivant $Q$, i.e., on consid\`ere les tours de Rokhlin
de m\^eme hauteur $H$ et de base $B\cap q$ pour chaque $q\in\bigvee_{k=0}^{H-1}T^{-k}Q$ tel que
l'intersection soit de mesure non-nulle. 

Ces sous-tours de Rokhlin sont $Q$-pures par construction. Elles portent donc chacune un $Q$-nom bien d\'efini. En
appliquant (C2), apr\`es \'elimination des sous-tours atypiques, on voit que ce $Q$-nom est essentiellement une
juxtaposition de (longues) $\db$-copies de d\'ebuts $P^{[m,n]}(B\cap q)$ d'une m\^eme suite $\omega$. Ces intervalles
$[n,m]$ d\'efinissent des tours de Rokhlin (de base $T^m(B\cap q)$ et de hauteur $n-m$). Celles-ci peuvent \^etre
juxtapos\'ees en une tour irr\'eguli\`ere essentiellement $Q$-pure qui satisfait \`a la condition~(C3).
 \qed

\subheading{(C3) $\implies$ (C4)}
C'est l'\'etape la plus importante (et la plus d\'elicate). Notre d\'emonstration est un raffinement de la
preuve donn\'ee dans le cas du rang un \cite{\KINa}.

On suppose que $(X,\Cal B,T,\mu)$ v\'erifie (C3). Soit $Q_1,Q_2,\dots$ une suite de
partitions de plus en plus fines dont l'union engendre la tribu des mesurables (ce qu'on notera:
$Q_n\nearrow\Cal B$).

\proclaim{Etape 1} 
Construction d'une suite de tours embo\^{\i}t\'ees $\tau_1,\tau_2,\dots$ telles que $P_n\eqdef P_{\tau_n}$
$\eps_n$-raffine $Q_n$ avec $\eps_n\searrow 0$.
\endproclaim

(on convient d'abr\'eger $P_{\tau_n}, h_{\tau'_n},\dots$ en $P_n,h'_n,\dots$).

Il suffira d'appliquer le lemme suivant avec, \`a la $n$i\`eme \'etape, $\delta=\delta_n$ un terme de s\'erie sommable:

 \proclaim{\lem}
Soit $\tau_1,\tau_2,\dots,\tau_n$ une suite de $n\geq0$ tours telles que $P_k$ $\eps_k$-raffine $Q_k$ avec $\eps_k>0$,
pour $k=1,\dots,n$. Soit $\delta>0$.

Il existe alors une suite de $n+1$ tours embo\^{\i}t\'ees $\tau'_1,\dots,\tau'_{n+1}$ telles que $P'_{n+1}$
$\delta$-raffine $Q_{n+1}$ et:
 $$
     d(P'_k,P_k) \eqdef 2\sum_i \mu(N'_k(i)\diffsym N_k(i)) < \delta \qquad \text{pour tout }1\leq k\leq n.
 $$
En particulier, $P'_k$ $(\eps_k+\delta)$-raffine $Q_k$.
 \endproclaim

En effet, le lemme fournit une suite de suites de tours irr\'eguli\`eres $\tau_k^{(n)}$, $1\leq k\leq
n$, $n=1,2,\dots$ embo\^{\i}t\'ees et telles que, pour $k$ fix\'e, $\tau_k^{(n)}$ converge quand $n\to\infty$ dans un
sens \'evident vers une tour $\tau_k$. La suite limite est bien \'evidemment embo\^{\i}t\'ee et $P_k$ $\sum_{l\geq
k}\delta_l$-raffine $Q_k$. La preuve du lemme ach\`evera donc l'\'etape~1.

\medbreak

Prouvons le lemme. L'id\'ee est la suivante: on prend une tour irr\'eguli\`ere $\tau'_{n+1}$ dont la partition soit si
fine que la base de $\tau_n$ est essentiellement une union de niveaux de $\tau'_{n+1}$, i.e., $\tau_n$ et $\tau'_{n+1}$
sont d\'ej\`a presque embo\^{\i}t\'ees: il suffira d'une correction tr\`es petite de $\tau_n$ en $\tau'_n$ pour obtenir
un embo\^{\i}tement exact. Enfin, pour maintenir les embo\^{\i}tements des $\tau_1,\dots,\tau_n$ il suffira de d\'efinir
$\tau'_k$ en fonction de $\tau'_{k+1}$ avec la m\^eme formule que $\tau_k$ en fonction de $\tau_{k+1}$. 

Voyons les d\'etails. Soit $H=\max h_n$ et $H_1 = 8 \delta^{-1} H^2$. 

\proclaim{\clm} 
On peut trouver une tour $\tau'_{n+1}$ dont la partition associ\'ee, $P'_{n+1}$,
$\delta/8H$-raffine la  partition $\Pi=\{h_n^{-1}(j):j=0,\dots,H-1\}\vee Q_{n+1}$ et telle que:
 $$
       \Delta'_{n+1} \eqdef \min_{1\leq s\leq S'_{n+1}} {k'}^{n+1}_s-{k'}^{n+1}_{s-1} \geq H_1
 $$
o\`u ${k'}^n_1<\dots<{k'}^n_{S'_{n+1}}$ sont les hauteurs de $\tau'_{n+1}$ (on convient que
${k'}^{n+1}_0=0$).
\endproclaim

\demo{Preuve de l'affirmation}
Posons $H_2=48 \delta^{-1} H\cdot H_1$ et $\gamma=\delta/48H H_2^{2}$. Soit
$R$ une partition dont les \'el\'ements sont de mesure au plus $\gamma/2$. D'apr\`es~(C3), il existe
une tour $\tau'$ dont la partition associ\'ee $\gamma/2$-raffine $\Pi\vee R$.

La mesure de chaque niveau de $\tau'$ est au plus \'egale \`a $\gamma$. On peut donc
supprimer les sous-tours de hauteur inf\'erieure \`a $H_2$: il y en a au plus $H_2$ de mesure chacune
au plus $H_2\cdot \gamma$. Ce faisant, on supprime une mesure au plus \'egale \`a $(H_2)^2\gamma=
\delta/48H$. 

Pour avoir la minoration sur $\Delta'_{n+1}$, il suffit maintenant de remplacer $h'$ par $H_1[h'/H_1]$ o\`u
$[\cdot]$ repr\'esente la partie enti\`ere. Ce faisant on supprime une mesure au plus \'egale \`a $H_1/\min
\{h'(x)>0:x\in X\} \leq H_1/H_2= \delta/48H$. Soit $\tau_{n+1}'$ la tour r\'esultant de ces suppressions.

La partition associ\'ee \`a $\tau'$ $\gamma/2$-raffine $\Pi$ avec $\gamma/2\leq \delta/24H$. La partition associ\'ee \`a
$\tau'_{n+1}$ raffine donc $\Pi$ \`a la pr\'ecision: $\delta/24H+2\mu(\text{enlev\'e})= \delta/8H$. 

L'affirmation est d\'emontr\'ee.
 \qed \enddemo

\demo{Preuve du lemme}
Soit donc $\tau'_{n+1}$ comme dans l'affirmation. Modifions $\tau_n$ en $\tau'_n$ embo\^{\i}t\'ee
dans $\tau'_{n+1}$. Pour cela on d\'efinit $\tau'_n$ de sorte que:
 \roster
  \item"(i) " pour chaque $j=1,\dots,\max h_n$, ${h'}_n^{-1}(j)$ est une union d'\'el\'ements de
$P'_{n+1}$.
  \item"(ii)" $\sum_j \mu({h'}_n^{-1}(j)\diffsym h_n^{-1}(j))<\delta/8H$.
 \endroster
Pour avoir un embo\^{\i}tement, on doit supprimer de l'union mentionn\'ee au point~(i) d\'efinissant ${h'}_n^{-1}(j)$
d'abord $X\setminus\tau'_{n+1}$ et ensuite les niveaux $N'_{n+1}(l)$ tels que:
 $$
     T^iN'_{n+1}(l)\ne N'_{n+1}(l+i) \text{ pour un }i\in[1,j[,
 $$
i.e., tels que ${h'}^{-1}_{n+1}([l+1,l+j[) \ne \emptyset$. Ce faisant on accro\^{\i}t la
distance~(ii) d'au plus $\delta/8H+H/H_1$. La distance~(ii) apr\`es cette correction est donc
major\'ee par $\delta/2H$.

On en d\'eduit que:
 $$
     d(P'_n,P_n)\leq 2 H \sum_j \mu({h'}_n^{-1}(j)\diffsym h_n^{-1}(j))
           \leq \delta.
 \tag\equ$$

$\tau'_n$ a donc les bonnes propri\'et\'es. Il reste \`a maintenir l'embo\^{\i}tement de $\tau_n$ dans $\tau_{n-1}$,
etc. Les tours $\tau_n$ et $\tau_{n-1}$ \'etaient d\'ej\`a embo\^{\i}t\'ees: la base de chaque sous-tour de Rokhlin de
$\tau_{n-1}$ \'etait une union de niveaux de $\tau_n$:
 $$
          h_{n-1}^{-1}(j) = \bigcup_{i=1}^{I_j} N_n(l_{ij}) \qquad \text{pour }j=1,2,\dots,\max h_{n-1}
 $$
chaque $l_{ij}$ v\'erifiant: $T^m(N_n(l_{ij}))=N_n(l_{ij}+m)$ pour $m\in[1,j[$, i.e.,
$h_n^{-1}([l_{ij}+1,l_{ij}+j[)=\emptyset$. 

Pour pr\'eserver cet embo\^{\i}tement, alors qu'on a d\^u modifier $\tau_n$ (mais sans rajouter de nouvelles hauteurs),
on garde la formule, en y substituant $\tau'_n$ \`a $\tau_n$. On obtient ainsi $\tau'_{n-1}$ avec $d(P'_{n-1},P_{n-1})
\leq d(P'_n,P_n)\leq \delta$. On r\'epercute ensuite cette nouvelle modification sur $\tau_{n-2}$,
etc.

Le lemme est d\'emontr\'e.\qed
 \enddemo

\proclaim{Etape 2}
On peut trouver une partition $P$ g\'en\'eratrice pour laquelle toutes les tours $\tau_n$
soient pures (quitte \`a remplacer les tours $\tau_n$ par des tours embo\^{\i}t\'ees $\tau'_n$ avec
$d(P_n',P_n)\to0$).  
\endproclaim

Pour garantir la puret\'e des tours, on construit chaque \'el\'ement de $P$ comme une union de
certains niveaux pris dans les tours $\tau_n$. Pour que $P$ soit g\'en\'eratrice, il suffit que,
pour $\mu$-presque tout $x$, la connaissance du nom $P^\Z(x)$ suffise \`a d\'eterminer le niveau dans lequel $x$ se
trouve par rapport \`a chaque tour $\tau_n$. En effet ceci d\'etermine $Q_n(x)$, l'\'el\'ement de $Q_n$ contenant $x$,
pour $x$ en dehors d'un ensemble de mesure au plus $\eps_n$: en supposant $\sum_n \eps_n<\infty$ et en appliquant
Borel-Cantelli, on voit que ceci d\'etermine $Q_n(x)$ pour $n$ assez grand. Mais $Q_n\nearrow\Cal B$. 

On souhaite poser $P=\{A,B,C,D\}$ avec:
 $$\align
        A &= \bigcup_{n\geq1} N_{n}(k_1^n-n-2)\cup N_{n}(k_1^n-1) \\
        B &= \bigcup_{n\geq1} \bigcup_{s\geq 2} N_{n}(k_s^n-n-2)\cup N_{n}(k_s^n-1)\\
        C &= \bigcup_{n\geq1} \bigcup_{s\geq 1} N_{n}(k_s^n-n-1)\cup\dots\cup N_{n}(k_s^n-2)\\
        D &= X\setminus(A\cup B\cup C)\\
 \endalign$$
o\`u $0<k_1^n<k_2^n<\dots<k_{S_n}^n=\max h_n$ sont les hauteurs de $\tau_n$, i.e., les entiers $k>0$ tels que
$h_n^{-1}(k)\ne\emptyset$. Si c'est possible, alors on lit sur le $P$-nom d'un point les instants o\`u il franchit la
base et le to\^{\i}t de chacune des tours $\tau_n$, et donc le niveau de $\tau_n$ contenant $x$ (ou si $x\in\tau_n$). 

En effet, il suffit de consid\'erer les apparitions des mots $AC^nA$ et $BC^nB$. Par exemple, l'apparition de $AC^nA\dots
BC^nB\dots BC^nB\dots AC^nA$ (les points de suspensions repr\'esentant des symboles $D$ ou des mots
$AC^mA$ ou $BC^mB$ avec $m\ne n$) signale le passage dans la sous-tour de $\tau_n$ de hauteur
$k_n^3$. Plus pr\'ecis\'ement, si on a:
 $$\matrix\format\l&\r&\r&\r&\r&\r&\r&\r&\r&\r\\
P^\Z(x):\qquad & \dots\quad &
     AC^nA & \quad\dots\quad &    BC^nB & \quad\dots\quad &    BC^nB &  \quad\dots\quad &    AC^nA &
\quad\dots \\
  &  &
  \uparrow &                 & \uparrow &                 & \uparrow &                  & \uparrow &
\\
  &  &
       t_1 &                 &      t_2 &                 &      t_3 &                 &       t_4 &
\\
 \endmatrix$$
cela correspond \`a une travers\'ee de la tour $\tau_n$ entre les instants $t_1-k_n^1+1$ (niveau $0$)
et $t_3$ (to\^{\i}t au niveau $k_n^3-1$).

Il reste \`a montrer, qu'apr\`es modification des tours, la d\'efinition propos\'ee pour $P$ est possible. Il suffit
que, pour chaque $n\geq1$:
 \roster
  \item"(i)" $k_1^n-n-2\geq H \eqdef \max h_{n-1}$.
  \item"(ii)" $\bigcup_{s=1}^{S_n}\bigcup_{m=1}^{n+2} N_n(k_n^s-m)$ est disjoint des tours de rang $<n$ (pour avoir la
disjonction de $A,B,C$).
 \endroster
La condition~(i) est facilement satisfaite (cf.\ l'\'etape~1).

Pour satisfaire la condition~(ii), il suffit de retirer des tours de rang $<n$ les niveaux en question, ainsi que leurs
pr\'e-images par $T,\dots,T^{H-1}$. Vu~(i), ces pr\'e-images ne sont autres que les niveaux $N_n(k_n^s-m)$ pour
$m\in[1,n+2+H]$.

La mesure \`a retirer des tours d'ordre $<n$ repr\'esente donc une fraction de $\tau_n$ (donc une mesure) major\'ee par
$(n+2+H)/\Delta_n$, l'\'ecart minimal entre deux hauteurs de $\tau_n$. Mais $\Delta_n\nearrow\infty$ aussi vite que l'on
veut. On peut donc effectuer ces corrections, tout en pr\'eservant l'embo\^{\i}tement en remaniant les tours
pr\'ec\'edentes comme dans l'\'etape~1.

Ceci ach\`eve l'\'etape~2 et la preuve de (C3)$\implies$(C4).

\remark{\rem} \edef\remSupport{\therem}
1. On peut se contenter d'une partition $P$ \`a deux \'el\'ements (il suffit de ``coder'' les nombres de $0$ \`a $3$ en
base $2$).

2. On peut ins\'erer dans les tours irr\'eguli\`eres tous les mots finis sur $P$ que l'on veut, et donc faire
appara\^\i tre ces mots dans la suite pistante.
\endrem\endremark

\subheading{(C4) $\implies$ $1$-pistage}
Supposons que $(X,\Cal B,T,\mu)$ v\'erifie (C4): on a donc une suite de tours $\tau_1,\tau_2,\dots$
qui sont embo\^{\i}t\'ees, pures pour une certaine partition g\'en\'eratrice $P$ et dont les partitions associ\'ees
v\'erifient: $P_n\eqdef P_{\tau_n}\nearrow\Cal B$.

Soit $\omega^{(n)}$ le $P$-nom de la tour $\tau_n$. Comme $P_n\nearrow\Cal B$ et que
l'espace est non-atomique, la mesure de $\tau_n$ tend vers $1$ tandis que celle de la base,
$N_n(0)$ tend vers z\'ero.  On en d\'eduit qu'il existe des entiers $l_n\to\infty$ tels que:
 $$
   \mu(B_n) \to 1 \qquad \text{avec }
       B_n\eqdef \bigcup_{s=1}^{S_n} \bigcup_{j=l_n}^{k_n^s-l_n} T^j( h_n^{-1}(k_n^s) )
 $$
En effet, $\mu(\tau_n\setminus B_n) = 2l_n\sum_{s=1}^{S_n} \mu( h_n^{-1}(k_n^s) ) = 2l_n
\mu(N_n(0))$, or $\mu(N_n(0))\to0$. Remarquons que $x\in B_n$ implique que
$P^{[-r,s]}(x)=\omega^{(n)}$ avec $r,s\geq l_n$.

En passant \`a une sous-suite, on peut supposer $\sum_{n\geq1} \mu(X\setminus B_n)<\infty$.
$\mu$-presque tout $x\in X$ est donc chaque $B_n$, d\`es que $n$ est assez grand: il existe
$r,s,n\to\infty$ tels que:
 $$
         P^{[-r,s]}(x) = \omega^{(n)}_0\dots\omega^{(n)}_{r+s}
 $$

On aura donc la propri\'et\'e de pistage si on peut faire en sorte que les $P$-nom des $\tau_n$, $\omega^{(n)}$, soient
les d\'ebuts d'une m\^eme suite $\omega\in P^\N$, i.e., qu'ils se prolongent les uns les autres. Vu les hypoth\`eses de
puret\'e, il suffit d'obtenir que les bases des tours $\tau_n$ et $\tau_{n+1}$ aient une intersection de mesure
non-nulle.

Remarquons que, comme dans la preuve de l'implication (C3)$\implies$(C4), on peut faire en sorte que
$\Delta_n$, la diff\'erence minimale entre deux hauteurs de $\tau_n$, tende vers l'infini. 

\proclaim{\lem}
Soit $\tau_1,\tau_2,\dots$ une suite de tours irr\'eguli\`eres $P$-pures, embo\^{\i}t\'ees et telles que $P_n\nearrow\Cal
B$ et $\Delta_n\nearrow\infty$. 

Fixons $\delta>0$. Pour $m$ assez grand, il existe $\tau'_n$ et $\tau'_m$ telles que:
 \roster
  \item $\tau'_n,\tau'_m$ sont encore $P$-pures et $N'_m(0)\subset N'_n(0)\subset N_n(0)$.
  \item $d(P'_n,P_n)<\delta$, $d(P'_m,P_m)<\delta$.
  \item $\tau'_m$ s'embo\^{\i}te dans $\tau'_n$.
 \endroster
En particulier, (1) implique que ${\omega'}^{(m)}$ prolonge ${\omega'}^{(n)}=\omega^{(n)}$.
\endproclaim

On en d\'eduira l'existence d'une suite de tours comme annonc\'e ((C4)+$P$-noms se prolongeant les uns les autres) en
proc\'edant comme dans la preuve de (C3)$\implies$(C4), apr\`es avoir remarqu\'e que la r\'epercution du d\'eplacement de
$\tau_n$ sur $\tau_1,\dots,\tau_{n-1}$ pr\'eserve la propri\'et\'e: $\omega^{(i)}$ prolonge
$\omega^{(j)}$ pour $1\leq j\leq i<n$.

\demo{Preuve}
On fixe $\delta>0$. En utilisant l'ergodicit\'e, on voit que:
 $$
     \mu\left(\bigcup_{k=0}^K T^{-k}N_n(0)\right)>1-\delta/6.
 $$ 
pour $K$ grand. Fixons $m>n$ suffisamment grand (on verra \`a quel point au cours de la
d\'emonstration).

L'union $U_m=\bigcup_{s=1}^{S_m} \bigcup_{i=1}^{K/\delta} N_m(k_m^s-i)$ est de mesure au plus
$K/\delta\Delta_m$ o\`u $\Delta_m$ est la plus petite diff\'erence entre deux hauteurs de $\tau_m$.
$m$ \'etant grand, $\Delta_m$ l'est aussi et on peut supposer que $U_m$ est de mesure au
plus $\delta/6$.
 
Soit $L$ minimal tel que $N_m(L)\cap\left(\bigcup_{k=0}^K T^{-k}N_n(0)\setminus U_m\right)$ soit de
mesure non-nulle. La mesure de l'union des niveaux de $\tau_m$ inf\'erieurs \`a $L$ est major\'ee par $(2/6)\delta$.
Enlevons-les et supposons d\'esormais $L=0$. La tour restante a une hauteur minimale au moins
\'egale \`a $K/\delta$, en effet, $U_m$ \'etant une union de niveaux, on a: $N_m(L)\subset X\setminus U_m$, et donc, par
d\'efinition de $U_m$, la hauteur restante est bien d'au moins $K/\delta$.

Soit $k\in[0,K[$ tel que $T^k(N_m(0))$ rencontre $N_n(0)$ sur un ensemble de mesure positive. 
Remarquons qu'on peut supposer que $K$ est plus petit que la hauteur minimale de $\tau_m$, d'o\`u:
$T^k(N_m(0))=N_m(k)$. Comme $\tau_m$ s'embo\^{\i}te dans $\tau_n$, on a en fait $N_m(k)\subset N_n(0)$.

Supprimons les niveaux en dessous de $k$. Ce faisant on enl\`eve une fraction de la mesure major\'ee par
$K/\Delta_m$ qu'on peut supposer major\'e par $\delta/6$. La tour $\tau'_m$ ainsi obtenue d\'efinit bien une suite finie
${\omega'}^{(m)}$ prolongeant $\omega^{(n)}$. D'autre part, $d(P'_m,P_m)<\delta$.

Comme dans la preuve de (C3)$\implies$(C4), on doit corriger $\tau_n$ en $\tau'_n$ pour maintenir
l'em\-bo\^{\i}\-te\-ment. Enfin les bases de $\tau_n$ et de $\tau'_n$ s'intersectent sur un ensemble de
mesure positive: les suites $\omega^{(n)}$ et ${\omega'}^{(n)}$ sont donc bien les m\^emes.
 \qed
 \enddemo

Le $1$-pistage est ainsi d\'emontr\'e, achevant la preuve du th\'eor\`eme~\thmCaract.
\qed

\head \newsec Cons\'equences \endhead

\subheading{\newsub Condition n\'ecessaire pour le $1$-pistage}

La caract\'erisation (C1) du th\'eor\`eme~\thmCaract (ou plut\^ot l'implication $1$-pistage$\implies$(C1)) donne
aussit\^ot une condition n\'ecessaire (un peu brutale!) pour le $1$-pistage:

\proclaim{\cor}
Soit $(X,\Cal B,T,\mu)$ un syst\`eme dynamique $1$-pist\'e. Soit $P$ une partition quelconque.
Pour $w=w_0\dots w_{n-1}\in P^n$, on note $B_\db(w,\eps)=\{x\in X: \db(P^{[0,n[}(x),w)<\eps\}$.
On a alors, pour tout $\eps>0$ et tout $N<\infty$,
 $$
    \bigcup_{n\geq N} T^{n-1}B_\db(\omega_0\dots\omega_{n-1},\eps) = X  \mod{\mu}, 
 $$
pour une certaine suite $\omega\in P^\N$.
En particulier:
 $$
       \sum_{n\geq1} \max_{w\in P^n} \mu(B_\db(w,\eps)) = \infty \qquad \forall \eps>0.
 $$
\endproclaim

\subheading{\newsub $1$-pistage et rang local}
D'apr\`es le th\'eor\`eme~\thmCaract, il est manifeste que rang $1$ implique $1$-pist\'e. On a en fait beaucoup
plus: la propri\'et\'e de rang local \cite{\FERa}, qui g\'en\'eralise strictement le rang fini, suffit \`a
entra\^{\i}ner le $1$-pistage.

Rappelons que le rang local peut \^etre caract\'eris\'e \cite{\KINa} par l'existence de tours $\theta_1,\theta_2,\dots$
 \roster
  \item"(i)" deux-\`a-deux embo\^{\i}t\'ees.
  \item"(ii)" de hauteurs $h_n$ tendant vers l'infini. 
  \item"(iii)" de mesure minor\'ee par une constante $a>0$.
  \item"(iv)" pures par rapport \`a  une partition g\'en\'eratrice fix\'ee.
 \endroster

\proclaim{\prop}
Le rang local implique le $1$-pistage.
\endproclaim

On utilise:

\proclaim{\lem}  \edef\lemRgLoc{\thelem}
Si $(X,\Cal B,T,\mu)$ est de rang local, alors on peut supposer que les tours de Rokhlin $\theta_1,\theta_2,\dots$
ci-dessus v\'erifient de plus:
 \roster
  \item pour chaque $k\geq1$,
     $${\mu(\theta_{n+1} \setminus (\theta_k\cup\dots\cup\theta_n))}
               / {\mu(\theta_{n+1})\mu(X \setminus \theta_k\cup\dots\cup\theta_n)))}
       \overset{n\to\infty}\to{\longrightarrow} 1$$
 (on convient de ce que $0/0=1$).
  \item $h_{n+1}/h_{n}\to \infty$.
 \endroster
\endproclaim

\demo{Preuve du lemme}
On utilise la caract\'erisation du rang local rappel\'ee ci-dessus. En passant \`a une sous-suite (ii)
donne (2). La seule chose \`a v\'erifier est le point (1). Soit $n\geq10$. Voyons qu'on peut supposer que l'\'ecart dans
(1) par rapport \`a la limite est major\'e par $1/n$ pour $1\leq k\leq n$ (pour une sous-suite l\'eg\`erement
modifi\'ee).

Soit $U_k=\theta_k\cup\dots\cup\theta_n$. Si $\mu(U_k)=1$, il n'y a rien \`a d\'emontrer. Supposons donc $\mu(U_k)<1$.
D'apr\`es le th\'eor\`eme ergodique de Birkhoff, il existe $Y$ de mesure au moins $1-a/n$ et $L<\infty$ tel que si $x\in
Y$ et $m\geq L$ alors, pour chaque $1\leq k\leq n$:
 $$
     \frac1m\sum_{p=0}^{m-1} \chi_{X\setminus U_k}(T^px) = (1+1/n)^{\pm1} \mu(X\setminus U_k)>0
 \tag\equ$$\edef\equfreq{\theequ}
($x^{\pm1}$ repr\'ensentant un nombre entre $x^{-1}$ et $x^{+1}$).

Quitte \`a sauter un certain nombre de tours dans la suite $\theta_{n+1},\theta_{n+2},\dots$ on peut supposer que
$h_{n+1}\geq 2 L$. Quitte \`a enlever de $\theta_{n+1}$ une fraction major\'ee en mesure et en hauteur par $1/n$, on
peut supposer que son premier niveau rencontre $Y$ sur un ensemble de mesure positive.  La hauteur restante est
minor\'ee par $L$, donc, vu l'embo\^{\i}tement, \equfreq\ donne, pour chaque $1\leq k\leq n$:
 $$
     \left| \frac {\mu(\theta_{n+1} \setminus (\theta_k\cup\dots\cup\theta_n))}
            {\mu(\theta_{n+1})\mu(X \setminus \theta_k\cup\dots\cup\theta_n)))} - 1 \right|
        \leq \frac 1n.
 $$
On peut donc d\'efinir la sous-suite de tours par r\'ecurrence sur $n$.
 \qed\enddemo

\demo{Preuve de la proposition}
Soit $T$ de rang local. On applique le lemme pr\'ec\'edent. Notons $P$ la partition g\'en\'eratrice pour laquelle les
tours sont pures.  Soit $\omega^n\in P^{h_n}$ d\'efini par: le $i$\`eme niveau de $\theta_n$ est inclus dans
$\omega^n_i$. Soit $\omega\in P^\N$ d\'efini par $\omega_{h_n}\dots\omega_{h_{n+1}-1}=\omega^{n+1}_{h_n}\dots
\omega^{n+1}_{h_{n+1}-1}$. 

Montrons que $\mu(\bigcup_{n\geq N} \theta_n)=1$ quel que soit $N<\infty$. Posons $U_m=\bigcup_{k=N}^m \theta_k$ et
$\eps_m=|1-\mu(\theta_{m+1}\setminus U_m)/\mu(\theta_{m+1})\mu(X\setminus U_m)|$. On calcule:
 $$\align
    \mu\left(X\setminus U_{m+1} \right) 
          & = 1 - \left( \mu(U_m) + \mu(\theta_{m+1}\setminus U_m) \right) \\
          & \leq 1 - \mu(U_m) - (1-\eps_m) \mu(\theta_{m+1})\mu(X\setminus U_m )\\
          & \leq \mu(X\setminus U_m )(1 - (1-\eps_m) \mu(\theta_{m+1}) ) 
 \endalign$$
Comme $\sum_m \mu(\theta_m)=\infty$, on voit que $\mu(U_m)$ tend vers $1$. On en d\'eduit: $\mu$-presque tout $x\in
X$ est dans $\theta_m$ pour des $m$ arbitrairement grands. Autrement dit: $P^{[-k,-k+h_m]}(x)=\omega^m$. Donc:
$P^{[-k+h_{m-1},-k+h_m]}(x)=\omega_{h_{m-1}}\dots\omega_{h_m-1}$ et:
 $$
       \db(P^{[-k,-k+h_m]}(x),\omega_0\dots\omega_{h_m}) \leq 2h_{m-1}/h_m \to 0.
 $$
Ci-dessus, $k\geq0$, $-k+h_m\geq0$ et $h_m\to\infty$. On a donc le $1$-pistage d'apr\`es la caract\'erisation~(C1) du
th\'eor\`eme~\thmCaract.
\qed
\enddemo

\subheading{\newsub Description topologique}

On voudrait munir l'ensemble des syst\`emes d'entropie nulle d'une notion to\-po\-lo\-gi\-que d'en\-sem\-ble
n\'egligeable. Les ensembles de premi\`ere cat\'egorie c'est-\`a-dire union d'une collection d\'enombrable de ferm\'es
d'int\'erieur vide fournissent une telle notion pourvu que l'espace lui-m\^eme ne soit pas de premi\`ere cat\'egorie.
D'apr\`es le th\'eor\`eme de Baire, il suffit que l'espace soit un espace m\'etrique complet.

Pour avoir une telle structure, on consid\`ere plut\^ot que les syst\`emes dynamiques les processus stochastiques.
Rappelons qu'un \new{processus stochastique} sur $N$ symboles est simplement la donn\'ee d'une mesure de probabilit\'e
invariante du d\'ecalage sur $N$ symboles. Tout syst\`eme dynamique probabiliste $T$ muni d'une partition $\alpha$
(finie et mesurable) ordonn\'ee $P=\{P_0,\dots,P_{N-1}\}$ d\'efinit un tel processus not\'e $(T,P)$.

On dira que le processus $(T,P)$ est $1$-pist\'e si le syst\`eme dynamique induit sur le d\'ecalage est $1$-pist\'e
comme syst\`eme dynamique probabiliste. Soulignons qu'on ne demande pas qu'il soit $1$-pist\'e par rapport \`a la
partition canonique du d\'ecalage.
 
Une "bonne" distance est la distance $\db$ fournie par la distance de Hamming entre suites finies.
Soit $(X,T,\mu,\alpha)$ et $(Y,S,\nu,\beta)$ deux processus. Si $\card\alpha\ne\card\beta$, alors
$\db((T,\alpha),(S,\beta))=1$. Sinon:
 $$
     \db((T,\alpha),(S,\beta)) = \sup_{n=1,2,\dots} \inf_\phi \frac1n \sum_{k=0}^{n-1} d(T^{-k}\alpha,\phi S^{-k}\beta)
 $$
o\`u:
 \roster
  \item $\phi:Y\to X$ d\'ecrit les isomorphismes entre espaces de probabilit\'e (ne pr\'eservant pas n\'ecessairement
la dynamique).
  \item $d(\alpha,\alpha')=\sum_{i=1}^N \mu(\alpha_i\diffsym\alpha'_i)$.
 \endroster

Remarquons que $\db$ n'induit pas une distance entre syst\`emes dynamiques probabilistes. En particulier si $T$ et $S$
sont deux syst\`emes d'entropie nulle alors il existe des partitions g\'en\'eratrices $\alpha,\beta$ \`a deux
\'el\'ements telles que $\db((T,\alpha),(S,\beta))$ est arbitrairement petite.

\proclaim{\prop\ (par exemple \cite{\PETa,\STa})} \edef\propDb{\theprop}
$\db((T,\alpha),(S,\beta))=0$ ssi les processus sont identiques. L'ensemble des processus ergodiques d'entropie nulle
est un espace m\'etrique complet pour $\db$.

$\db((T,\alpha),(S,\beta))<\eps$ est \'equivalente \`a n'importe laquelle des deux conditions suivantes :
 \roster
  \item le nom-$T,\alpha$ de $\mu$-presque tout point peut \^etre transform\'e en un nom g\'e\-n\'e\-ri\-que pour
$S,\beta$ en ne modifiant ce nom que sur un ensemble de fr\'equence $<\eps$.
  \item il existe un nom g\'en\'erique pour $T,\alpha$ et un nom g\'en\'erique pour $S,\beta$ qui ne diff\`erent que sur
une sous-suite de fr\'equence $<\eps$.
 \endroster
\endproclaim

On veut montrer:

\proclaim{\prop}
L'ensemble des syst\`emes $1$-pist\'es est un ensemble topologiquement n\'e\-gli\-gea\-ble et m\^eme un ferm\'e
d'int\'erieur vide dans la classe des processus ergodiques d'entropie nulle munie de $\db$.
\endproclaim

On admet l'existence d'un syst\`eme d'entropie nulle non $1$-pist\'e, syst\`eme dont nous construirons un exemple dans
la section~5.

\demo{Preuve}
Tout d'abord, montrons que l'ensemble des processus $1$-pist\'es est ferm\'e. 

Soit $(T_1,\alpha_1),(T_2,\alpha_2),\dots$ tendant vers $(T,\alpha)$ au sens de $\db$. Supposons $\alpha$ et les
$\alpha_n$ g\'en\'eratrices et chaque $T_n$ $1$-pist\'e et montrons que $T$ l'est alors aussi. On utilise la
caract\'erisation (C2) (th\'eor\`eme~\thmCaract). Soit $0<\eps<1$. Soit $n$ tel que $\db(T_n,T)<\eps^2/4$.

$T_n$ \'etant pist\'e,  il existe $\omega\in(\alpha_n)^\N$ tel que le nom $\alpha_n,T_n$ de $\mu_n$-presque tout point
peut \^etre $(1-\eps/2)$-recouvert par des $\db-\eps/2$-copies de d\'ebuts de $\omega$ de longueur au moins
$6\eps^{-1}N$. D'apr\`es la proposition~\propDb, on en d\'eduit que le nom $\alpha,T$ de $\mu$-presque tout point peut
\^etre $(1-\eps)$-recouvert par des $\db-\eps$-copies de d\'ebuts de $\omega$, identifi\'ee \`a une suite
\`a valeurs dans $\alpha$, de longueur au moins $2\eps^{-1}N$.

En utilisant \`a nouveau (C2), on voit que le syst\`eme $T$, donc le processus $(T,P)$, est bien $1$-pist\'e.
L'ensemble des processus $1$-pist\'es est bien ferm\'e. Pour voir qu'il est d'in\-t\'e\-rieur vide, on admet
provisoirement l'existence d'un syst\`eme $(U,\rho)$ d'entropie nulle qui n'est pas $1$-pist\'e (la section~\secNonPiste\
en  donne un exemple).

Consid\'erons un processus ergodique $1$-pist\'e $(T,\mu,P)$. On suppose $P$ g\'en\'eratrice. Soit $\eps>0$. Il s'agit de
construire un processus $(S,Q)$ avec $\db((T,P),(S,Q))<\eps$ et $S$ non $1$-pist\'e.

Tout d'abord, quitte \`a remplacer $(T,P)$ par un processus arbitrairement proche au sens de $\db$, on peut supposer que,
pour un entier $r$ assez grand, il y a au plus $\card P^r/(\card Q+1)$ cylindres de longueur $r$ de $(T,P)$ ayant
une mesure non-nulle. On note $P_r$ l'ensemble de ces cylindres. 

En effet, $T$ \'etant d'entropie nulle, le cardinal minimal d'une collection de $r$-cylindres dont l'union a une mesure
sup\'erieure \`a $1-\eps$ est plus petit que $e^{\eps r}<\card P^r/r(\card Q+1)$ pour $r$ assez grand. Fixons un tel $r$
et une collection $P_r$ correspondante. Consid\'erons la projection $\pi:P^\Z\to P^\Z$ induite par une application par
bloc $P^r\to P_r$ co\"\i ncidant avec l'identit\'e sur $P_r$. Posons $\nu=\frac1r\sum_{s=0}^{r-1} \sigma^s \pi(\mu)$ (on
a identifi\'e $X$ \`a $P^\Z$). Le processus $(P^\Z,\sigma,\nu)$ a au plus $r\card P_r$ $r$-cylindres de mesure
non-nulle: il r\'epond au probl\`eme.

On part de:
 $$\matrix\format\r&\r&\quad\c\quad&\l\\
     T'\;:\; & X\times\{0,\dots,M-1\}\times Y & \to & X\times\{0,\dots,M-1\}\times Y\\
             & (x,k,y) & \longmapsto & (Tx,k+1\mod M,U^{\delta_{0k}}(y))
 \endmatrix$$
avec $M=r\lceil \eps^{-1} \rceil$ et $\delta_{ij}=1$ si $i=j$, $0$ sinon. On munit $T'$ de la mesure $\mu\times n\times
\rho$ avec $n$ la mesure de comptage normalis\'ee. 

Fixons une injection $i:P_r\times (Q\cup\{*\})\to P^r$. On consid\`ere l'application:
 $$
     (x,k,y)\mapsto \left\{ \matrix\format\l\qquad &\l\\
                                   i(P^r(x)\times *)    & \text{si }k\ne0 \\
                                   i(P^r(x)\times Q(y)) & \text{si }k=0. \\
          \endmatrix\right.
 $$
Elle induit un codage, non $T'$-invariant,  $\psi$ de $X\times\{0,\dots,M-1\}\times Y$ dans $(P_r)^\Z\subset P^\Z$. On
pose $\mu'=\frac1r\sum_{s=0}^{r-1} \sigma^s \psi(\mu\times n\times\rho)$. Soit $S$ le d\'ecalage muni de
$\mu'$ et $Q$ la partition canonique de $P^\Z$. Clairement, $(S,Q)$ est ergodique et d'entropie nulle et v\'erifie:
$\db((X,P),(S,Q))<\eps$.

Supposons $S$ $1$-pist\'e et tirons-en une contradiction.

Tout d'abord, tout it\'er\'e $S^k$ serait encore pist\'e. En effet, si $S$ est $(Q',\Sigma)$-pist\'e alors $S^k$
est $({Q'}^k,\Sigma^k)$-pist\'e avec:
 $$
     \Sigma^k=\{(\omega_i\omega_{i+1}\dots\omega_{i+k-1})(\omega_{i+k}\dots)\dots\in ({Q'}^k)^\N:i=0,\dots,k-1\}.
 $$
$\Sigma^k$ \'etant un ensemble fini, c'est dire que $S^k$ est encore $1$-pist\'e. Donc $S^M$ est $1$-pist\'e.

Mais $U$ est un facteur de l'it\'er\'e $S^M$. Voyons que ceci implique que $U$ est $1$-pist\'e, donc une contradiction.
En effet, utilisons la caract\'erisation (C1). Soit $R$ une partition g\'en\'eratrice $U$. Il
suffit d'appliquer (C1) \`a la partition $\pi^{-1}U$ de $S^M$, puis de projeter pour obtenir (C1) pour
$U$ par rapport \`a $R$. $S$ n'est donc pas $1$-pist\'e. \qed
\enddemo

\subheading{\newsub Preuve du th\'eor\`eme~\thmLien}

\demo{Suites pistantes et suites g\'en\'eriques}
Pour voir qu'une suite pistante $\omega$ est quasi-g\'en\'erique, il suffit d'appliquer le th\'eor\`eme ergodique \`a $T$
et $T^{-1}$: on obtient que, pour presque tout point, d\`es que les entiers $n,m$ sont assez grands, les
fr\'equences le long du segment $[-n,m]$ de l'orbite sont bonnes. En comparant avec la d\'efinition du $1$-pistage, on
voit qu'il existe une suite d'entiers, de la forme $n_i+m_i$ avec les notations de cette
d\'efinition, tels que les fr\'equences de $\omega|[0,n_i+m_i]$ tendent vers les bonnes fr\'equences. Autrement dit, la
suite pistante est quasi-g\'en\'erique.

Pour voir que $\omega$ n'est pas forc\'ement g\'en\'erique, il suffit de consid\'erer les suites construites ci-dessous
pistant simultan\'ement des mesures distinctes. Ces suites, quasi-g\'en\'eriques pour plus d'une mesure, sont
n\'e\-ces\-sai\-re\-ment non-g\'en\'eriques. \qed
\enddemo

\demo{Abondance topologique}
Soit $G_\ell$ l'ensemble des suites de $\{0,1\}^\N$ dont le d\'ebut est de la forme $BXB^{m}$ avec $B$ une suite de
$\ell$ symboles, $X$ une suite finie quelconque sur $\{0,1\}$ et $m$ un entier suffisamment grand pour que
$(m+1)\ell>\ell^2|X|$ et $m\geq\ell^2$. 

$G_\ell$ est manifestement un ouvert dense de l'espace de Baire $\{0,1\}^\N$. Donc
$G=\bigcap_{\ell\geq1} G_\ell$ est un $G_\delta$-dense.

Soit $\omega\in G$. Consid\'erons le syst\`eme de rang un d\'efini symboliquement par
$B_0=\omega_0$ et $B_{n+1}=B_nX_nB_n^{m}$ avec $X_n$ d\'efini par:
 $$
     \omega_0\dots\omega_{k-1} =B_nX_nB_n^m \qquad\text{avec } k=(m+1)|B_n|+|X_n|,\; m\geq|B_n|^2
 $$ 
est l'\'ecriture la plus courte manifestant l'appartenance de $\omega$ \`a $G_\ell$ avec
$\ell=|B_n|\geq n$. Ceci d\'efinit bien un syst\`eme de rang un en raison de l'in\'egalit\'e: 
$|X_n|\leq n^{-2}|B_{n+1}|$. 

Tout segment fini extrait du nom de $\mu$-presque tout point co\"{\i}ncide avec un extrait de $B_n$ pour tout $n$
assez grand. Mais chaque $B_n$ est un d\'ebut de $\omega$. Donc ce syst\`eme de rang un est bien pist\'e par $\omega$
pour la partition canonique. \qed
\enddemo

\demo{Abondance au sens de la mesure}
Soit $T$ de rang local. Soit $P$ une partition g\'en\'eratrice par rapport \`a laquelle $T$ soit pist\'e.
Le rang local peut \^etre d\'efini par le fait qu'une fraction uniform\'ement minor\'ee du $P$-nom de presque
tout point consiste en des copies exactes d'un mot arbitrairement long. En proc\'edant comme dans la preuve de
du th\'eor\`eme~\thmCaract\ et en utilisant la remarque~\remPurete, on d\'eduit de la propri\'et\'e de rang
local l'existence d'une suite infinie de tours de Rokhlin qui sont:
 \roster
  \item"(i)" deux-\`a-deux embo\^{\i}t\'ees.
  \item"(ii)" de hauteurs tendant vers l'infini.
  \item"(iii)" de mesures minor\'ees par une constante $a>0$.
  \item"(iv)" chacune pure par rapport \`a $P$, la partition donn\'ee ci-dessus.
 \endroster
Enfin, on proc\`ede comme pour le lemme~\lemRgLoc\ pour voir qu'on peut aussi supposer:
 \roster
  \item $\Cal A_n$ d\'esignant l'alg\`ebre finie engendr\'ee par les partitions associ\'ees \`a $\theta_1,\dots,
\theta_n$:
   $$
          \sup_{A\in\Cal A_n} \left| \frac{\mu(A\cap\theta_{n+1})}{\mu(A)\mu(\theta_{n+1})} - 1 \right| 
              \overset{n\to\infty}\to{\longrightarrow} 0.
   $$
  \item $h_{n+1}/h_{n}\to \infty$.
 \endroster
La propri\'et\'e d'ind\'ependance (1) permet de voir que $\mu$-presque tout $x$ appartient \`a $\bigcup_{k=0}^{h_n/n}
T^k B_n$ pour une infinit\'e de $n$, si $B_n$, $h_n$ sont la base et la hauteur de $\theta_n$, en utilisant $\prod_n
(1-a/n) =0$.

$x$ \'etant fix\'e, on peut donc supposer qu'il appartient  \`a $\bigcup_{k=0}^{h_n/n} T^k B_n$ pour {\bf tout}
$n\geq1$, en passant \`a une sous-suite.

Posons $\theta_n'=\bigcup_{k=2h_n/n}^{h_n} T^k B_n$. $y\in \theta'_n$ implique que $P^{[-a,b]}(y)=P^{[0,a+b]}(x)$ avec
$a\geq h_n/n$ et $b\geq0$. La propri\'et\'e (1) donne comme ci-dessus que $\bigcup_{n\geq N} \theta'_n= X$ modulo $\mu$,
pour tout $N$ (en effet $\mu(\theta_n')/\mu(\theta_n)\to1$).  Donc $\mu$-presque tout $y$ est dans une infinit\'e de
$\theta'_n$. On a donc le $1$-pistage par le $P$-nom de $x$ d'apr\`es la caract\'erisation (C1).
\qed\enddemo

\demo{Co-pistage}
On consid\`ere donc deux syst\`emes dynamiques $1$-pist\'es ap\'eriodiques et ergodiques $T$ et $T'$. On veut trouver
un entier $N$, une partition g\'en\'eratrice \`a $N$ \'el\'ements pour chacun ($P$ et $P'$) et une suite
$\alpha\in\{0,1,\dots,N-1\}^\N$ telles que les syst\`emes soient pist\'es par rapport \`a
$(P,P_{\alpha_0}P_{\alpha_1}\dots)$, resp.\ $(P',P'_{\alpha_0}P'_{\alpha_1}\dots)$.

\medbreak

Comme soulign\'e dans la remarque~\remSupport, la construction d'une suite de tours embo\^\i t\'ees permet d'obtenir une
partition g\'en\'eratrice \`a deux \'el\'ements qui est admissible par le pistage et qui induit un codage envoyant la
mesure sur une mesure ayant le d\'ecalage tout entier pour support. On peut supposer que les deux syst\`emes sont
port\'es par le d\'ecalage $\{0,1,\dots,N-1\}^\Z$, ont pour m\^eme support (l'espace tout entier) et sont $1$-pist\'es
par rapport \`a la partition canonique.

Remarquons tout d'abord que le cylindre d\'efini par un mot quelconque $a_0\dots a_n$ est de $\mu$-mesure positive ssi
$a_0\dots a_n$ apparait dans $\omega$. L'\'egalit\'e des supports entra\^{\i}ne donc que les m\^emes mots finis
apparaissent dans $\omega$ et dans $\omega'$.

Construisons maintenant une suite $\alpha$ pistant les deux syst\`emes.

Supposons que pour un $n\geq1$, on a d\'ej\`a trouv\'e:
 $$
     \alpha_0\dots\alpha_{L-1}=\omega_d\dots\omega_{d+L-1}=\omega'_{d'}\dots\omega'_{d'+L-1}
            \qquad d,d'\geq0 \text{ et } L\geq n-1.
 $$
tel qu'on puisse recouvrir une fraction au moins $1-1/n$ du nom de  $\mu$- et $\nu$-presque tout point par des
copies exactes de segments initiaux de longueur $\geq n-1$ de $\alpha_0\dots\alpha_{L-1}$. Remarquons que cette assertion
est triviale pour $n=1$.

Soit un entier $\ell$ suffisamment grand pour que des copies exactes de d\'ebuts de $\omega$ de
longueurs comprises entre $\ell_-\eqdef 2(n+1)d$ et $\ell$ permettent de recouvrir $1-1/2(n+1)$ de $\mu$-presque toute
orbite. Prolongeons $\alpha$ en posant:
 $$
     \alpha_0\dots\alpha_{K-1} = \omega_{d}\dots\omega_{\ell-1} \qquad \text{avec } K=\ell-d.
 $$
On peut ainsi recouvrir $1-1/2(n+1)-d/\ell_-=1/(n+1)$ de $\mu$-presque toute orbite avec des d\'ebuts de
$\alpha_0\dots\alpha_{K-1}$ de longueur au moins $n$.

D'apr\`es la remarque ci-dessus, $\alpha_0\dots\alpha_{K-1}$ apparaissant dans $\omega$ appara\^{\i}t aussi dans
$\omega'$: il existe $d''$ tel que $\omega'_{d''}\dots\omega'_{d''+K-1}$ co\"{\i}ncide avec ce mot.

On proc\`ede alors comme ci-dessus en rempla\c{c}ant $d+L$ par $d''+K$, $\mu$ par $\nu$ et $\omega$ par $\omega'$.
On obtient un prolongement $\alpha_0\dots\alpha_{L'-1}$ avec les propri\'et\'es voulues au rang $n$.

On obtient donc par r\'ecurrence une suite infinie $\alpha$ qui piste \`a la fois $\mu$ et $\nu$.
\qed
\enddemo

\remark{\rem}
La m\^eme d\'emonstration se g\'en\'eralise \`a une infinit\'e d\'enombrable de syst\`emes $1$-pist\'es
ergodiques et ap\'eriodiques. Par contre, le cas d'une infinit\'e non-d\'enombrable ou encore d'une suite universelle
reste entier.
\endrem\endremark

\head \newsec Un syst\`eme $1$-pist\'e non-L.B. \endhead

L'exemple construit dans cette section est en partie inspir\'ee par la construction par J.~Feldman d'un syst\`eme
d'entropie nulle non-L.B.\ \cite{\FELa}.
\medbreak

Rappelons la d\'efinition de la propri\'et\'e L.B. (cf.\ \cite{\ORWa}) (on se restreint dans cet article aux syst\`emes
d'entropie nulle). C'est la fermeture, par induction et suspension, de l'ensemble des rotations irrationnelles. Le point
de vue le plus commode pour nous est fourni par la caract\'erisation symbolique donn\'ee par J.~Feldman
\cite{\FELa} en terme de la \new{distance $\fb$} entre suites finies: 
 $$
         \fb(a,b) = \frac{\vareps}{|a|+|b|}
 $$
o\`u $\vareps$ est le nombre minimal de symboles \`a effacer dans $a$ et dans $b$ pour obtenir des suites identiques,
i.e., c'est le plus petit entier $\vareps$ tel qu'il existe deux suites de longueur $r=(|a|+|b|-\vareps)/2$: $0\leq
i_1<\dots<i_r<|a|$ et $0\leq j_1<\dots<j_r<|b|$, telles que $a_{i_1}\dots a_{i_r}=b_{i_1}\dots b_{i_r}$. On dit que
$(i,j)$ est un \new{couplage} de \new{longueur} $r$.

\proclaim{\defn\ (J.~Feldman)}
Soit $(X,\Cal B,T,\mu)$. Si $P$ est une partition de $X$ et $n\geq1$, on d\'efinit la pseudo-distance $\fb_P^n(x,y)=
\fb(P^{[0,n[}(x),P^{[0,n[}(y))$.

$T$ est L.B.\ ssi pour toute partition $P$ et tout $\eps>0$, pour tout $n$ assez grand,
il existe une boule pour la pseudo-distance $\fb_P^n$ de rayon $\eps$ dont la mesure soit $\geq 1-\eps$.
\endproclaim

On montre \cite{\ORWa} qu'il est \'equivalent de demander la propri\'et\'e ci-dessus pour {\sl une} partition
g\'en\'eratrice.

\subheading{\newsub Construction par blocs}
La d\'emarche est standard. On la rappelle pour fixer les notations:
\medbreak

Soit $\Cal A$ un ensemble fini contenant le symbole $0$. Pour chaque ordre $n\geq n_0$, on d\'efinit une collection de
\new{types de blocs} $\{B_n^1,\dots,B_n^{N(n)}\}$. Chaque $B_n^i$ est une suite finie $B_n^i(0)\dots
B_n^i(l_n^i-1)$:
 \roster
  \item"---" \`a valeurs dans $\Cal A\setminus\{0\}$ si $n=n_0$.
  \item"---" \`a valeurs dans $\Cal B_{n-1}=\{0,B_{n-1}^1,\dots,B_{n-1}^{N(n-1)}\}$ si $n>n_0$.
 \endroster
On pose $\Cal B_{n_0-1}=\Cal A$.

Si $B$ est une suite de termes pris dans $\Cal B_n$, on note $\phi(B)$ la suite sur $\Cal B_{n-1}$ obtenue en
juxtaposant les termes de $B$. Pour $m\leq n$, on appelle $m$-d\'eveloppement de $B$ la suite $\phi^{n-m}(B)$ sur $\Cal
B_m$. Enfin le terme de \new{d\'eveloppement}, sans autre pr\'ecision, d\'esigne le $n_0$-d\'eveloppement. C'est une
suite sur $\Cal A$. 

On note $b_n^i$ le d\'eveloppement de $B_n^i$.

\medbreak

Notons $p(B|B')$ la fraction de la longueur occup\'ee par le d\'eveloppement de $B$ dans le d\'eveloppement
de $B'$.

Supposons que chaque $B_n^i$ appara\^{\i}t au moins une fois dans chaque $B_{n+1}^j$ et que:
 $$
    \sum_{n\geq n_0} \sup_k p(0|B_n^k) < \infty \text{ et }
    \sup_{n\geq n_0} \sup_{i,k,l}\frac{p(B_n^i|B_{n+1}^k)}{p(B_n^i|B_{n+1}^l)} < \infty.
 $$

Il existe alors un unique automorphisme d'espace de Lebesgue $(X,\Cal B,T,\mu)$ tel que, \`a chaque type de bloc
$B_n^i$ est associ\'ee une tour de Rokhlin de base $R_n^i$ et de hauteur $|b_n^i|$, avec les propri\'et\'es
suivantes:
 \roster
  \item l'union des tours d'ordre $n$, $\bigcup_{i=1}^{N(n)}\bigcup_{k=0}^{|b_n^i|-1} T^kR_n^i$ est de mesure tendant
vers $1$ quand $n\to\infty$.
  \item embo{\^\i}tement: $T^\ell R_n^j$ est l'union des $T^k R_{n+1}^i$ tels que le $k$i\`eme symbole du
d\'eveloppement de $B_{n+1}^i$ correspond au $\ell$i\`eme symbole du d\'eveloppement de $B_n^j$.
 \endroster
On en d\'eduit que:
 $$
         \mu(R^j_n)|b_n^j| = \lim_{m\to\infty} p(B_n^j|B_m^{i_m})
 $$
quelle que soit la suite des entiers $i_m\in\{1,\dots,N(m)\}$.
De plus, $T$ est ergodique.

La partition standard d\'efinie par une telle construction est $\{P_A:A\in\Cal A\}$ avec: $P_A$ l'union des
$T^kR_{n_0}^i$ si le $k$i\`eme symbole de $B_{n_0}^i$ est $A\in\Cal A\setminus\{0\}$, et $P_0$ le reste de
l'espace $X$. 

\subheading{\newsub D\'efinition de l'exemple}
Pour motiver cette d\'efinition, faisons quelques remarques.

$\bullet$ Pour garantir le $1$-pistage d'un syst\`eme $(\Cal A^\Z,\sigma,\mu)$ d\'efini comme ci-dessus, il suffit
d'avoir, pour tout $n\geq n_0$:
 \roster
  \item $p(B_n^1|B_{n+1}^i)\geq 1/n$ pour tout $i$. 
  \item le type $B_{n+1}^1$ commence par $B_n^1$ pour tout $n\geq n_0$.
  \item $\min_i |b_n^i|\to\infty$.
 \endroster

\demo{Preuve du $1$-pistage}
Gr\^ace \`a (2), on peut d\'efinir $\omega\in\Cal A^\N$ par:
 $$
      \omega|[0,|b_n^1|[ = b_n^1 \qquad \text{pour chaque }n\geq n_0.
 $$

Comme $\prod_n (1-1/n)=0$, on d\'eduit de (1) que, pour $\mu$-presque tout $x\in\Cal A^\Z$, il existe
$n_1<n_2<\dots<n_p<\dots$ tels que $x$ est dans la tour associ\'ee \`a $B_{n_p}^1$, en particulier:
$P^{[-a_p,b_p]}(x)=b_{n_p}^1$ avec $a_p,b_p\geq0$ et $a_p+b_p=|b_{n_p}^l|\to\infty$.

La caract\'erisation~(C1) du th\'eor\`eme~\thmCaract\ implique que $(\sigma,\mu)$ est $\omega$-pist\'e par rapport
\`a la partition $P$. Le syst\`eme obtenu en quotientant $(X,\Cal B,T,\mu)$ par la tribu $\bigvee_{k\in\Z} T^kP$ est
donc $1$-pist\'e. 

\qed \enddemo

$\bullet$ Pour emp\^echer la propri\'et\'e L.B.\ pour ce m\^eme syst\`eme quotient il suffit de maintenir l'\'ecart
entre les blocs de types diff\'erents, i.e., il suffit (cf. \cite{\ROTa, lemma~2.6}) de trouver une
constante $\delta_*>0$ tel que:
 $$
       \inf_{i\ne j} \fb(b_n^i,b_n^j) \geq \delta_* > 0 \qquad \text{pour tout $n\geq n_0$}.
 $$

Pour obtenir une telle minoration, on proc\`ede bien s\^ur par r\'ecurrence. Mais un couplage entre deux $n+1$-blocs
n'est pas forc\'ement compatible avec le d\'ecoupage de chacun des deux $n+1$-blocs en leurs $n$-blocs. On sera donc
amen\'e \`a minorer:
 $$
     \fb(n) \eqdef \inf_{u,v} \fb(u,v)
 $$
o\`u $u$ et $v$ parcourent les segments de longueurs $\geq |b_n|/n^2$ extraits de $n$-blocs de types distincts,
$|b_n|$ \'etant la longueur commune des $b_n^i$, quel que soit $i$.

Tr\`es grossi\`erement, si on consid\`ere un couplage entre $u$ et $v$ comme ci-dessus, alors la force de ce couplage est
donn\'ee par:
 $$
     [1-\fb(u,v)] \leq (1-\gamma)\left\{ (1-\delta) [1-\fb(n)] + \delta \cdot 1 \right\}
             = (1-\gamma)\left\{ 1-(1-\delta)\fb(n) \right\}
 $$
avec $\gamma$ la proportion des $n-1$-blocs que le couplage consid\'er\'e efface int\'egralement et $\delta$ la
proportion de $n-1$-blocs de m\^eme type qu'il met en correspondance.

Pour que montrer que $\fb(n)$ ne tend pas vers z\'ero, la meilleure majoration possible sur $\delta$ seul, $\delta\leq
\max_i p(B_n^i) = 1/n$, est insuffisante. Mais si les occurrences d'un type de $n$-bloc donn\'e dans $b_{n+1}^i$ et
celles dans $b_{n+1}^j$ sont suffisamment enchev\^etr\'ees, alors on ne pourra  en coupler une fraction significative
($\delta\geq n^{-3/2}$) sans avoir aussi \`a effacer beaucoup de $n$-blocs ($\gamma\geq n^{-3/2}$). 

Cet enchev\^etrement sera obtenu en disposant les $B_n^t$ de fa\c{c}on {\bf ind\'ependante} dans $B_{n+1}^i$ et
$B_{n+1}^j$ si $i\ne j$. 

\medbreak

On obtient ainsi la d\'efinition suivante pour notre syst\`eme. A l'ordre $n+1$, on d\'efinit $N(n+1)\eqdef n!+1$
types de blocs en posant:
 $$
       B_{n+1}^i = \left(C_{n+1}^i \right)^{\beta_n^{N(n+1)-i+1}}  \text{ pour }i=1,\dots,N(n+1).
 $$
avec:
 $$
       C_{n+1}^i = \underbrace{B_n^1\dots\dots\dots B_n^1}_{\beta_n^{i-1}\nu_n}\;
                             \underbrace{B_n^2\dots B_n^2}_{\beta_n^{i-1}}\; \dots \;
                             \underbrace{B_n^{N(n)}\dots B_n^{N(n)}}_{\beta_n^{i-1}}
 $$
o\`u $\beta_n=|C_{n+1}^1|=\nu_n+N(n)-1$ et $\nu_n=(N(n)-1)/(n-1)$.

En particulier, le $k$i\`eme bloc dans $B_{n+1}^i$ est de type $t\ne1$ (resp.\ de type $1$) si $t-1+\nu_n$ (resp.\ %
si un symbole $<\nu_n$) est le $i$i\`eme chiffre de $k-1$ en base $\beta_n$: ceci fournit l'ind\'ependance recherch\'ee.

Remarquons que le choix de $\nu_n$ garantit que $p(B_n^1|B_{n+1}^i)=1/n$ et que le nombre de
r\'ep\'etitions de $C_{n+1}^i$ est toujours tr\`es grand ($\geq\beta$) et que la longueur de $B_{n+1}^i$ est
ind\'ependante de $i$. On note: $|B_n|\eqdef |B_n^i|$ pour tout $i$. Remarquons que:
 $$\gather
    |C_{n+1}^k| = \beta_n^{k-1} (N(n)-1+\nu_n) \\
    |B_{n+1}|=|b_{n+1}|/|b_n| = \beta_n^{N(n+1)+1}=\left((N(n)-1)\cdot\frac{n}{n-1}\right)^{N(n+1)+1}. \\
 \endgather$$
Cette quantit\'e cro\^{\i}t en fonction de $n$ de fa\c{c}on super-exponentielle.

\subheading{\newsub Preuve du caract\`ere non-L.B}
Elle repose sur les deux lemmes suivants.

Le premier permet de traiter le cas des couplages entre $n+1$-blocs de type distincts dont une proportion significative
des appariements se produisent entre $n$-blocs de m\^eme type. C'est ici qu'on utilise la structure particuli\`ere de
notre exemple et en particulier l'enchev\^etrement.

\proclaim{\lem\ (enchev\^etrement)} \edef\lemEnchev{\thelem}
Soit $(I,J)$ un couplage entre $U=C_{n+1}^k$ et $V=\left( C_{n+1}^{k'} \right)^q$ avec $k>k'$ et $q\geq1$. Notons $M$
la longueur de ce couplage. 

Supposons que $\kappa\eqdef \frac{2M}{|U|+|V|}\geq200n^2/\beta_n$. Alors les indices des appariements entre
symboles $B_n^1$ dans $U$ et dans $V$ forment un intervalle, i.e., $(U\circ I)^{-1}(B_n^1)=[M',M'']$ avec:
 $$\gather
        M'=1 \text { et } M''\geq \left(1-\frac1{n^2}\right) M \\
        I(M')\leq \frac1n|U| \text{ et } J(M')\geq \frac n2\kappa|V|
 \endgather$$
\endproclaim

Le deuxi\`eme lemme traite le cas des autres couplages.

\proclaim{\lem} \edef\lemCouplage{\thelem}
Soit $U$ et $V$ des segments extraits de $B_{n+1}^k$ et $B_{n+1}^{k'}$. Soit $u,v$
les d\'eveloppements de $U$ et de $V$. Supposons qu'un des couplages $(i,j)$ r\'ealisant $\fb(u,v)$ ne fasse
correspondre que peu de symboles figurant dans des $n$-blocs de m\^eme type, i.e., supposons que:
 $$
       R = \card\{s:U([i_s/|b_n|])=V([j_s/|b_n|])\} \leq \frac{1}{n^{3/2}} \frac{|U|+|V|}{2}.
 $$
Alors:
 $$
       \fb(u,v) \geq \min\left( \fb(n) - \frac{\cte}{n^{3/2}}, \frac15 \right).
 $$
\endproclaim

Appliquons ces deux lemmes avant de les d\'emontrer. 

\proclaim{\clm}
Pour tout $n\geq n_0$,
 $$
     \fb(n+1) \geq \max\left(\fb(n)-\frac{\cte}{n^{3/2}},\frac15\right).
 $$
\endproclaim

\demo{Preuve de l'affirmation}
Soit $u,v$ comme dans la d\'efinition de $\fb(n+1)$: des segments de longueurs $\geq|b_{n+1}|/n^2$ extraits de
$n+1$-blocs $b_{n+1}^k, b_{n+1}^{k'}$ avec $k<k'$. Quitte \`a changer $\fb(u,v)$ d'au plus $|C_{n+1}^{N(n+1)}|
/(n^2|B_{n+1}|)= n^2 \beta_n^{-1}=o(n^{-3/2})$, on peut supposer que $u$ est le d\'eveloppement de $U$, la
concat\'enation d'un certain nombre de $C_{n+1}^k$.

Ecrivons $u=u^1\dots u^q$, o\`u chaque $u^p$ est une copie du d\'eveloppement de $C_{n+1}^k$. D\'ecoupons $v=v^1\dots
v^q$ de fa\c{c}on \new{compatible modulo $(i,j)$}, i.e., si $i_m$ tombe dans $u^p$, alors $j_m$ tombe dans $v^p$. On a
donc:
 $$
    \fb(u,v) = \sum_{p=1}^q \frac{|u^p|+|v^p|}{|u|+|v|} \fb(u^p,v^p).
 $$
Il suffit donc de minorer $\fb(u',v')$ avec $u'$ le d\'eveloppement de $C_{n+1}^k$ et $v'$ extrait du d\'eveloppement de
$C_{n+1}^{k'}$ r\'ep\'et\'e ind\'efiniment. Si $|v'|\leq(2/3)|u'|$, alors $\fb(u',v')\geq1/5$ et il n'y a rien d'autre
\`a d\'emontrer. On peut supposer que $|v'|\geq(2/3)|u'|$. La longueur $|v'|$ est donc au moins de l'ordre de $\beta_n
|C_{n+1}^{k'}|$. On peut donc supposer \'egalement que $v'$ est le d\'eveloppement de $V'=(C_{n+1}^{k'})^Q$ avec $Q$
tr\`es grand.

On omet d\'esormais les primes sur $u,v,U,V$.

Fixons $(i,j)$ un couplage entre $u$ et $v$ r\'ealisant la distance $\fb(u,v)$. Soit $(I,J)$ le couplage induit entre
$U$ et $V$, i.e., en notant $M$ sa longueur:
 $$
    \{(I_m,J_m):m=1,\dots,M\} = \{(i_s,j_s):U([i_s/|b_n|])=V([j_s/|b_n|])\}.
 $$

Si $2M/(|U|+|V|)<n^{-3/2}$, il suffit d'appliquer le lemme~\lemCouplage. Supposons maintenant: $2M/(|U|+|V|)\geq n^{-3/2}
\geq 200n^2/\beta_n$. \smallbreak

D\'ecoupons $U=U'U''$, $U'$ regroupant les symboles $B_n^1$ et $U''$ les autres. Ceci
induit un d\'ecoupage $u=u'u''$. Soit $V=V'V''$ tel que le d\'ecoupage induit $v=v'v''$ soit compatible modulo $(i,j)$
(ce qui est plus fort que modulo $(I,J)$). 

D'apr\`es le lemme~\lemEnchev, on a $|U'|=|U|/n$ et $|V'|\geq \frac n2 \frac{2}{n^{3/2}} |U| = |U|/n^{1/2}$. Donc
$\fb(u',v')\geq 1-\cte/n^{1/2}\geq 1/2$ et:
 $$
     \fb(u,v) \geq \frac{|u'|+|v'|}{|u|+|v|} \frac 12 + \frac{|u''|+|v''|}{|u|+|v|} \fb(u'',v'')
 $$
Maintenant, toujours d'apr\`es le lemme~\lemEnchev, le couplage entre $U''$ et $V''$ est de longueur au plus
$M/n^2\leq(|U|+|V|)/n^2$, toujours \`a cause du lemme~\lemEnchev. On peut donc appliquer le lemme~\lemCouplage, qui
donne la minoration voulue.

L'affirmation est d\'emontr\'ee. Le syst\`eme n'est pas L.B.
\qed
\enddemo

\demo{Preuve du lemme \lemEnchev}
On omet l'indice $n$ l\`a o\`u cela ne provoque pas de confusion: on \'ecrira donc $N,\beta,\nu,\dots$ au lieu de
$N(n),\beta_n,\nu_n,\dots$.

Remarquons que, par construction, l'application $t\mapsto U_i=C_{n+1}^k(i)$ est croissante (en identifiant $B_n^i$ \`a
l'entier $i$). On en d\'eduit que, pour chaque $i$, $(U\circ I)^{-1}(B_n^i)$ est un intervalle. Notons $a_i$ le nombre
d'appariements entre symboles $B_n^i$. Par construction: $0\leq a_1\leq \nu\beta^{k-1}$ et, pour $i=2,\dots,N$,
$0\leq a_i\leq \beta^{k-1}$. Soit:
 $$
      \Cal I \eqdef \{ i=1,\dots,N : a_i\geq 10 \nu(i) \beta^{k'-1} \}
 $$
avec $\nu(1)=\nu_n$ et $\nu(i)=1$ pour $i=2,\dots,N$. 

On calcule: $\sum_{i\notin\Cal I} a_i \leq 10\beta^{k'-1} (N-1+\nu) = 10 \beta^{k'-k}|U|$. Par hypoth\`ese, $k-k'\geq 1$,
donc:
 $$
      \sum_{i\notin\Cal I} a_i \leq 10\beta^{-1} |U|
            \leq \frac{1}{10n^2} \frac{100n^2}{\beta} \cdot 2 \frac{|U|+|V|}{2} \leq \frac{1}{10n^2} M.
 $$
D'o\`u:
 $$
      \sum_{i\in\Cal I} a_i = M - \sum_{i\notin\Cal I} a_i\geq \left(1-\frac{1}{10}\frac{1}{n^2}\right) M.
 $$
Soit $i\in\Cal I$. $(U\circ I)^{-1}(i)$ est un intervalle $[m_-,m_+]$ de $a_i$ entiers. Remarquons que les symboles
$B_n^i$ apparaissent dans $V$ par segments de longueur $\nu(i)\beta^{k'-1}$ se reproduisant \`a intervalles de longueur
$\beta^{k'}$. D'o\`u:
 $$
       J(m_+) - J(m_-) \geq \beta^{k'} \left( \left\lceil a_i/\nu(i)\beta^{k'-1}\right\rceil - 1 \right)
 $$
($\lceil x\rceil$ d\'esignant le plus petit entier $\geq x$). Mais $a_i/\nu(i)\beta^{k'-1}\geq 10$ (car $i\in\Cal I$)
donc:
 $$
       J(m_+) - J(m_-) \geq \frac{9}{10} \frac{\beta}{\nu(i)} a_i.
 \tag\equ$$ \edef\equJ{\theequ}
Si on avait:
 $$
      a_2+\dots+a_N\geq M/n^2
 \tag\equ$$ \edef\equAbsurde{\theequ}
alors on aurait: $\sum_{i\in\Cal I\setminus1} a_i\geq \sum_{i\ne1} a_i - \sum_{i\notin\Cal I} a_i \geq \frac{9}{10}
M/n^2$ et donc:
 $$\align
       J(M)-J(1) &\geq \sum_{i\in\Cal I\setminus1} J(m_+(i))-J(m_-(i))
                  \geq \frac{9}{10} \beta \sum_{i\in\Cal I\setminus1}  a_i\\
                 &\geq \frac{9}{10} \beta \cdot \frac{9}{10} \frac{M}{n^2} \geq 81 \cdot \frac{\beta M}{100 n^2} \\
                 &\geq 81 (|U|+|V|)>|V|,
 \endalign$$
ce qui est absurde. \equAbsurde\ est donc impossible. On a donc: $a_1>(1-1/n^2)M$. En particulier,
 $$
      a_1\geq \frac12 M \geq 50n^2\beta^{-1} |U| = 50 n^2 \beta^{-1} \cdot \beta^{k-1}(\nu+N-1)
         \geq 50 n^2 \nu \beta^{k-2} > 10 \nu \beta^{k'-1}
 $$
donc $1\in\Cal I$ et, vu \equJ:
 $$
     J(M'')-J(M') \geq \frac{9}{10} \frac{\beta}{\nu} \cdot a_1 
                  \geq \frac{9}{10} n \cdot (1-n^{-2}) M \geq \frac12 n \kappa |V|
 $$
($M=\frac12\kappa(|U|+|V|)$).

Enfin, $U^{-1}(B_n^1)=[0,\dots,|U|/n-1]$ d'o\`u $(U\circ I)^{-1}(B_n^1)=[M',M'']$ avec $M'=1$ et $I(M'')\leq |U|/n$.
\qed
\enddemo
 
\demo{Preuve du Lemme \lemCouplage}
Soit $u=u^1\dots u^s$ avec $u^t$ le d\'eveloppement de $U(t)$. Soit $v=v^1\dots v^s$ un d\'ecoupage de $v$ compatible,
modulo $(i,j)$, avec le d\'ecoupage pr\'ec\'edent de $u$. R\'eciproquement, la structure en $n$-blocs de $V$ induit un
d\'ecoupage sur $v$, et donc sur chaque $v^t$. On \'ecrit $v^t=v^t_1\dots v^t_q$, $q$ d\'ependant de $t$. Soit
$u^t=u^t_1\dots u^t_q$ un d\'ecoupage compatible. On a alors:
 $$
      \fb(u,v) = \sum_{t=1}^s \frac{|u^t|+|v^t|}{|u|+|v|} \fb(u^t,v^t)
               = \sum_{t=1}^s \sum_{p=1}^q \frac{|u^t_p|+|v^t_p|}{|u|+|v|} \fb(u^t_p,v^t_p)
 \tag\equ$$ \edef\equDist{\theequ}
On peut \'ecarter les termes mettant en jeu des {\bf suites de longueurs trop dif\-f\'e\-ren\-tes} qui donneront de
grandes distances ($\geq1/5$): si on pose $E_1 = \{ (t,p) : \frac{|u^t|}{|v^t|} \text{ ou }\frac{|u^t_p|}{|v^t_p|} 
\notin[2/3,3/2]\}$, alors:
 $$
        \sum_{(t,p)\in E_1} \frac{|u^t_p|+|v^t_p|}{|u|+|v|} \fb(u^t_p,u^t_q)
                    \geq \frac15 \frac{\sum_{(t,p)\in E_1} |u^t_p|+|v^t_p|}{|u|+|v|}.
 $$
Vu la minoration souhait\'ee, on peut se placer dans le pire des cas et supposer $E_1=\emptyset$.

On peut \'ecarter parmi les termes restants ceux dont les {\bf longueurs sont trop petites}, ils n'interviennent dans
\equDist\ que pour une fraction n\'egligeable: soit $E_2=\{ (t,p)\notin E_1 : \min(|u^t_p|,|v^t_p|)<|b_n|/n^2 \}$, alors:
 $$
     \sum_{(t,p)\in E_2} |u^t_p|+|v^t_p| < \card E_2 \cdot \frac52 \frac{|b_n|}{n^2}
           \leq 3|U| \cdot \frac52 \frac{|b_n|}{n^2} \leq 3\frac{|U|}{|B_n|} \cdot \frac52 \frac{|b_n|}{n^2}
           \leq 8 \frac{|u|+|v|}{n^2}
 $$
o\`u on a utilis\'e que si $(t,p)\notin E_1$ alors $v^t$ est de longueur au plus $(3/2) |b_n|$ et donc est \`a cheval sur
$q\leq3$ $n$-blocs. Les termes correspondants \`a $E_2$ n'interviennent dans \equDist\ que pour
$o(n^{-3/2})$.

Il reste donc les termes $\fb(u^t_p,v^t_p)$ avec $\frac{|u^t_p|}{|v^t_p|}\in[2/3,3/2]$ et $\min(|u^t_p|,|v^t_p|)\geq
|b_n|/n^2$. On les subdivise en deux cat\'egories: ceux tels que $u^t_p,v^t_p$ sont extraits d'un m\^eme type de
$n$-blocs et les autres. On note $F_1$, resp.\ $F_2$ l'ensemble des indices $(t,p)$ de la premi\`ere, resp.\ la seconde
cat\'egorie de termes.

Les premiers peuvent \^etre \`a distance nulle. Les autres sont, par d\'efinition, \`a distance au moins $\fb(n)$.
Finalement, \`a des termes en $o(n^{-3/2})$ pr\`es:
 $$
     \fb(u,v) \geq \fb(n) \left( 1 - \sum_{ (t,p)\in F_1 } \frac{|u^t_p|+|v^t_p|}{|u|+|v|} \right).
 $$
Mais, par hypoth\`ese, 
 $$
     \sum_{ (t,p)\in F_1 } \frac{|u^t_p|+|v^t_p|}{|u|+|v|} \leq \frac{2R}{|U|+|V|} \leq \frac{1}{n^{3/2}}.
 $$
Le lemme est d\'emontr\'e.\qed
\enddemo

\head \newsec Un syst\`eme L.B. d'entropie nulle non $1$-pist\'e \endhead

Les boules-$\db$ par rapport \`a une partition $P$ sont les ensembles:
 $$
    B_P(w_0\dots w_{n-1},\eps) = \{x\in X: \db(P^{[0,n[}(x),w_0\dots w_{n-1})<\eps\} \qquad w\in P^n.
 $$

\proclaim{\thm}
Pour tout $\Gamma>1$, il existe $(X,\Cal B,T,\mu)$ un automorphisme ergodique d'un espace de
Lebesgue qui est d'entropie nulle, l\^achement Bernoulli et qui admet une partition g\'en\'eratrice
$P$ telle que, pour $\eps>0$ assez petit, on a, pour tout $n$ assez grand:
 $$
     \mu(B_P(w_0\dots w_{n-1},\eps)) \leq n^{-\Gamma} \qquad \forall w\in P^n.
 $$
\endproclaim

\proclaim{\cor}
Il existe un syst\`eme d'entropie nulle, l\^achement Bernoulli qui n'est pas $1$-pist\'e.
\endproclaim

Le corollaire d\'ecoule du th\'eor\`eme et du crit\`ere \corCritere.

\remark{\rem}
On peut interpr\'eter ce r\'esultat en disant qu'\`a partir d'une rotation, pour laquelle l'exposant $\Gamma$
comme ci-dessus le plus grand possible est z\'ero, on peut, en induisant convenablement (cf.\ la d\'efinition des
syst\`emes L.B.\ au d\'ebut de la section pr\'ec\'edente), obtenir un exposant $\Gamma$ arbitrairement grand.
\endremark

\remark{\rem}
Pour une partition g\'en\'eratrice quelconque $Q$, le th\'eor\`eme admet le corollaire suivant. Pour tout $\eps>0$, il
existe $Y\subset X$ de mesure $1-\eps$ tel que, pour tout $n$ assez grand:
 $$
       \mu(B_Q(w_0\dots w_{n-1},\eps)\cap Y) \leq n^{-\Gamma} \qquad \forall w\in Q^n.
 $$
On observe que la complexit\'e mesur\'ee \cite{\FERb} d'un syst\`eme l\^achement Bernoulli d'en\-tro\-pie nulle peut
\^etre d'un ordre polynomial arbitrairement \'elev\'e.
 \endrem
 \endremark

On d\'efinit notre exemple par une construction par blocs comme dans la section pr\'ec\'edente.

Soit $n_0$ et $l_0$ deux grands entiers et $\gamma>\max(\Gamma,1)$. A chaque rang $n\geq n_0$, on se donne une
collection de blocs $B_n^1,\dots,B_n^{N_n}$ \'equiprobables. 

Pour $n=n_0$ ces blocs sont les suites d'entiers disjointes suivantes:
 $$
    B_{n_0}^i=[(i-1)l_0+1][(i-2)l_0+2]\dots[il_0] \qquad i=1,2,\dots,N_0\eqdef l_0^{\gamma-1}.
 $$

Pour $n+1>n_0$, on construit $N_{n+1}$ $n+1$-blocs de m\^eme longueur $l_{n+1}$ par concat\'enation de $n$-blocs et de
symboles $0$ d'apr\`es:
 $$
    B_{n+1}^i = B_n^1 0^{a_n^i(1)} B_n^2 0^{a_n^i(2)} \dots B_n^{N_n} 0^{a_n^i(N_n)}
    \qquad i=1,\dots,N_{n+1}.
 $$
Le nombre de $n+1$-blocs, $N_{n+1}$, sera pris \'egal \`a $l_{n+1}^{\gamma-1}$ pour obtenir
l'estimation recherch\'ee. En effet, les tours associ\'ees \`a chacun des $B_{n+1}^i$ seront de mesure $1/N_{n+1}$ et
donc chacun des niveaux sera de la mesure voulue, i.e., $1/N_{n+1}l_{n+1}=1/l_{n+1}^\gamma$.

On note $P$ la partition standard de $X$ en $N_0l_0+1$ sous-ensembles correspondant aux symboles $1,2,\dots,N_0l_0$ and
$0$ ci-dessus.

Si la connaissance de $P^\ell(x)$ \`a $\eps-\db$-pr\`es suffit \`a d\'eterminer la tour d'ordre $n+1$ et le
niveau contenant $x$ d\`es que $\ell\geq l_n^{1+\beta}$, alors on obtiendra une majoration de la mesure des boules-$\db$
en $l_{n+1}^{-\gamma}$ pour $l_n^{1+\beta}\leq\ell\leq l_{n+1}^{1+\beta}$, d'o\`u une majoration en
$\ell^{-\gamma/(1+\beta)}\leq \ell^{-\Gamma}$ pour $\beta>0$ assez petit.

Remarquons que l'entropie du syst\`eme ainsi d\'efini est nulle, vu $N_{n+1}=l_{n+1}^{\gamma-1}$: le nombre de blocs
cro\^\it polyn\^omialement avec la longueur.

Les nombres $a_n^i(k)$ seront choisis de mani\`ere \`a ce que:
 \roster
  \item"(C1)" tous les $n+1$-blocs soient de m\^eme longueur, $l_{n+1}=|b_{n+1}^i|$, pour tout $i$;
  \item"(C2)" les $0$ intercal\'es \`a l'ordre $n+1$ repr\'esentent au plus $1/n^2$ de cette
longueur. 
 \endroster

\medbreak

On obtient $(X,\Cal B,T,\mu)$ un automorphisme ergodique d'un espace de Lebesgue. 

Remarquons que la structure choisie implique que cet automorphisme est l\^a\-che\-ment Bernoulli, car en effa\c{c}ant les
``0'' intercalaires d'ordre $n$, (une modification inf\'erieure \`a $n^{-2}\to0$ pour la distance
$\fb$), on voit que tous les types de $n+1$-blocs se confondent, ce qui implique la propri\'et\'e L.B. d'apr\`es
\cite{\ROTa}.

\medbreak
{\sl Notation.}
$(B_n^i)_a^b$ est le mot fini sur $P$ extrait du d\'eveloppement $b_n^i$ \`a partir de la position $a$ (incluse) et
jusqu'\`a la position $b$ (excluse).

\subheading{\newsub Choix des intercalaires}
Rappelons que $\gamma>\max(\Gamma,1)$ et que $\beta>0$ est petit.

   \proclaim{\clm} \edef\clmSystNonPiste{\theclm}
On peut choisir les longueurs intercalaires $\{a_n^i(k)\}_{nik}$ de sorte qu'il y ait beaucoup de
blocs, i.e.:
 $$
    N_n \eqdef {l_n}^{\gamma-1}
 $$
et qu'il existe $\eps_*>0$ avec la propri\'et\'e suivante:
 \block
Pour tous $n,m\geq n_0$, si $\db((B_n^i)|_a^b,(B_m^j)|_c^d)<\eps_*$ et $b-a=d-c\geq
l'_n$ avec:
 $$
    l'_n \eqdef l_{n-1}^{1+\beta}
 $$ 
($l'_{n_0}=1$), alors $(B_m^j)_{c-a}^{c-a+l_n}$ est une occurrence explicite du bloc
$B_n^i$.
 \endblock
    \endproclaim

\remark{\rem}
D\'efinissons formellement la notion d'\new{occurrence explicite}. On proc\`ede par induction sur $m-n$.
On dit que $(B_m^j)_a^b$ est une occurrence explicite de $B_n^i$ si:
 \roster
  \item $n=m$ et $[a,b[=[0,l_m[$.
  \item $n=m-1$, $a=|b_n^10^{a_n^j(1)}\dots b_n^{i-1}0^{a_n^j(i-1)}|$ et $b=a+l_n$.
  \item $0\leq n<m-1$ et il existe des entiers $c,d,k$ tels que $c\leq a<b\leq d$ et 
$(B_m^j)_c^d$ est une occurrence explicite de $B_{m-1}^k$ et $(B_{m-1}^k)_{a-c}^{b-c}$ est
une occurrence explicite de $B_n^i$.
 \endroster
\endrem\endremark

On va disposer les $``0''$ intercalaires de fa\c{c}on \`a ce que:
 \roster
  \item"(C3)" la connaissance d'un nombre r\'eduit de co{\"\i}ncidences permette de lire le type du
$n$-bloc que l'on regarde ainsi que sa position exacte.
 \endroster

Remarquons que pour $n=n_0$ l'affirmation est satisfaite avec $\eps_*=1$.
Supposons fix\'ees les longueurs intercalaires $a_m^i(k)$ d'ordre $m<n$ de sorte qu'il existe
$\eps_n>0$ valable au rang $n$. On cherche des longueurs intercalaires $\{a_n^i(k)\}_{ik}$ et un
$\eps_{n+1}$ valable au rang $n+1$ qui ne soit pas trop petit, i.e., avec $\sum_n \log\eps_n/\eps_{n+1}<\infty$ (de
sorte qu'on puisse poser $\eps_*=\inf_n\eps_n>0$).

\subheading{1} 
La connaissance d'un extrait d'un $n+1$-bloc de longueur $l'_{n+1}$ doit d\'eterminer le
type et la position de ce $n+1$-bloc. On d\'ecoupe donc $B_{n+1}^i$ en \new{sous-blocs de base} comprenant $h_n$
$n$-blocs cons\'ecutifs avec:
 $$
      h_n \eqdef \frac{1}{n^2} \frac{l'_{n+1}}{l_n} 
        = \frac{1}{n^2} \cdot l_n^\beta
 $$
Tout extrait de longueur au moins $l_{n+1}'$ contient au moins $n^2$ sous-blocs de base.

Fixons $u=(B_{n+1}^i)_a^b$ et $v=(B_m^j)_c^d$ s\'epar\'es par:
 $$
    \db(u,v)<\eps_{n+1}\eqdef (1-3n^{-2})^3 \eps_n.
 $$
D\'ecoupons $u$ suivant les limites des sous-blocs de base de $B_{n+1}^i$. Notons $u=u_1\dots u_r$ le d\'ecoupage
obtenu. Soit $v=v_1\dots v_r$ le d\'ecoupage correspondant de $v$ (i.e., $|v_s|=|u_s|$ et on dit que
$v_s$ est le \new{vis-\`a-vis} de $u_s$). $u_1$ et $u_r$ sont peut-\^etre  des sous-blocs de base
incomplets. Comme $r\geq n^2$,
 $$
    \db(u,v) \geq (1-2n^{-2}) \frac1{r-2} \sum_{s=2}^{r-1} \db(u_s,v_s)
 $$
on peut trouver au moins un $s=2,\dots,r-1$ tel que $\db(u_s,v_s)<(1-3n^{-2})^2\eps_n$.

\subheading{2} 
On veut que chaque sous-bloc de base (et donc $u_s$) comporte toute l'in\-for\-ma\-tion. C'est pourquoi on prend
$a_n^i(\cdot)$ $h_n$-p\'eriodique. Pour caract\'eriser le $i$\`eme type de $n+1$-bloc, on se donnera ci-dessous une
fonction $f_n^i$. On code cette fonction par la longueur des intervalles intercalaires en posant:
 $$
    a_n^i(k+qh_n)=f_n^i(\lceil k/n^2\rceil)\in\{1,\dots,l_n^\beta\} \qquad \text{pour }k=1,\dots,h_n-1,\; q\geq 0
 \tag\equ$$ \edef\equAink{\theequ}
et  $a_n^i(qh_n)=h_n l_n^\beta - \sum_{k=1}^{h_n-1} a_n^i(k)$ pour $q\geq1$, pour avoir l'uniformit\'e des
longueurs des $n+1$-blocs. $\lceil x\rceil$ est le plus petit entier $\geq x$. On code donc $n^2$ fois chaque valeur de
$f_n^i$ dans chaque sous-bloc de base. On appelle chacune de ces suites de $n^2$ $n$-blocs ayant les m\^emes
intervalles intercalaires d'ordre $n$ un \new{segment r\'ep\'etitif}.

Evaluons $N_n$ et $l_n$. Vu $\sum_{k=1}^{h_n} a_n^i(k) = h_n l_n^\beta$, $l_{n+1}=N_nl_n(1+l_n^{-(1-\beta)})$. Or
$N_n=l_n^{\gamma-1}$ donc $l_{n+1}=l_n^{\gamma}(1+l_n^{-(1-\beta)})$. On en
d\'eduit:
 $$
      \log l_n = \cte\cdot \gamma^n (1+o(1))
 $$ 
avec $\cte$ une constante positive.

D'apr\`es \equAink, $u_s$ qui comprend $h_n$ $n$-blocs se d\'ecompose en $h_n/n^2$ segments r\'ep\'etitifs. Il
en existe au moins $h_n/n^4$ qui sont {\bf $\delta$-proches} de leurs vis-\`a-vis dans $v$: la distance-$\db$ entre les
deux est major\'ee par: $\delta\eqdef (1-n^{-2})^{-1}\db(u_s,v_s)\leq(1-3n^{-2})\eps_n$. 

\subheading{3}
Consid\'erons un des $h_n/n^4$ segments r\'ep\'etitifs de $u_s$ proches de leurs vis-\`a-vis. L'union (non contig\"ue) des $n^2$
$n$-blocs du segment (priv\'es des intervalles intercalaires) repr\'esente une fraction au moins
$(1+l_n^{-(1-\beta)})^{-1}\geq 1-n^{-2}$ de la longueur totale. La distance-$\db$ entre la restriction du segment \`a
cette union et le vis-\`a-vis dans $v$ est donc major\'ee par $(1-3n^{-2})/(1-n^{-2}) \cdot \eps_n$.

Enfin, en divisant l'union des $n^2$ $n$-blocs en chacun de ces $n$-blocs, on trouve au moins deux $n$-blocs, disons
$B_n^s$ et $B_n^r$, qui sont $\eps_n$-proches de leur vis-\`a-vis, vu $\frac{1-3n^{-2}}{(1-n^{-2})(1-2n^{-2})} <1$.

D'apr\`es l'hypoth\`ese de r\'ecurrence les vis-\`a-vis de ces deux $n$-blocs de $B_{n+1}^i$ sont des occurences
explicites de $n$-blocs du m\^eme type dans $B_m^j$. Notons $\Delta$ le nombre de symboles apparaissant entre les
d\'eveloppements de $B_n^s$ et de $B_n^r$ dans $u$. 

On a: $\Delta = (r-s-1)l_n + (r-s)f_n^i(\lceil s/n^2\rceil)$ en regardant $u$. On a la m\^eme formule avec $j$ en
regardant $v$. D'o\`u $f_n^i(\lceil s/n^2\rceil)=f_n^j(\lceil s/n^2\rceil)$. Les $h_n/n^4$ segments r\'ep\'etitifs
$\delta$-proches fournissent donc chacun une co\"\i ncidence distincte. Pour que ce soit suffisant pour
reconna\^{\i}tre le $n+1$-bloc, on demande que les fonctions $f_n^i$ satisfassent:
 $$
    i\ne j \implies \card\{1\leq k< h_n/n^2:f_n^i(k)=f_n^j(k)\} < h_n/n^4.
 \tag*$$
On aura alors la propri\'et\'e d'occurrence explicite annonc\'ee.

\subheading{\newsub Existence des $f_n^i$}
On v\'erifie qu'on peut bien trouver une collection de $N_{n+1}$ fonctions deux-\`a-deux
s\'epar\'ees au sens de (*). Posons:
 $$
   p=h_n/n^2-1,\quad b=l_n^\beta,\quad \delta=h_n/n^4.
 $$
On fixe $f_n^1:\{1,\dots,p\}\to\{1,\dots,b\}$ arbitrairement. Puis, $f_n^1,\dots,f_n^i$ \'etant
choisis, on veut pouvoir choisir $f_n^{i+1}$ distincte de toutes les fonctions $f:\{1,\dots,p\}\to\{1,\dots,b\}$
telles que:
 $$
    \card\{k:f(k)=f_n^j(k)\}\geq \delta
 $$
pour un $j=1,\dots,i$.

Notons $E$ le nombre maximal de fonctions ainsi exclues pour un $j$ donn\'e. On va montrer:
 $$
    (N_{n+1}-1) E < b^p
 $$
i.e., apr\`es avoir choisi $i<N_{n+1}$ fonctions, l'ensemble des fonctions exclues est encore strictement plus petit
que l'ensemble de toutes les fonctions: on peut donc continuer jusqu'\`a $N_{n+1}$.

En sommant sur le nombre $q$ de co\"\i ncidences, 
 $
   E = \sum_{q=\delta}^p C^q_p (b-1)^{p-q}.
 $
Or: 
 $$\multline
   \frac{C^{q+1}_p (b-1)^{p-q-1}}{C^q_p (b-1)^{p-q}}  = \frac{p-q}{(q+1)(b-1)} 
     <    \frac{p}{\delta(b-1)} \\
     =    \frac{h_n/n^2-1}{(h_n/n^4)(l_n^\beta-1)}
     \leq \cte \cdot \frac{n^2}{l_n^\beta}<1,
 \endmultline$$
vu l'estimation de $l_n$ ci-dessus.

Le premier terme de la somme $E$ majore donc les autres et on a:
 $$
    E \leq p \cdot C^\delta_p \cdot (b-1)^{p-\delta}
      \leq p \cdot \sqrt{\delta} \left(\frac {ep}\delta\right)^\delta \cdot b^{p-\delta}
 $$
en utilisant la formule de Stirling: $n!\sim (n/e)^n\sqrt{2\pi/n}$ et donc:
$C^\delta_p< p^\delta/\delta!< \sqrt{\delta}p^\delta/(\delta/e)^\delta=
 \sqrt{\delta}(ep/\delta)^\delta$. Finalement:
 $$\align
  \frac{b^p}{E} &\geq \frac1{p\sqrt{\delta}} \left(\frac{b \delta}{ep}\right)^\delta 
    \geq \frac1{h_n^{3/2}} \left( \frac{l_n^\beta}{3n^2} \right)^{h_n/n^4}
    \geq \frac{1}{l_n^{\frac32 \beta}}
       \left(\frac{l_n^{\beta}}{3n^2}\right)^{l_n^\beta/n^6}
\\
    &= \exp \left\{ \beta \frac{l_n^\beta}{n^6}\log l_n + \dots \right\}
 \endalign$$
(les termes n\'eglig\'es dans l'accolade sont des $o(\cdot)$ du terme \'ecrit). Or
 $$
    N_{n+1} \leq 2 N_nl_n = 2 l_n^{\gamma} = \exp \left\{ \gamma \log l_n + \dots \right\}.
 $$
On voit donc que $b^p/E>> N_{n+1}$.

\subheading{\newsub Majoration de $\mu(B_P(u,\eps))$}
Il reste \`a calculer la mesure des $\eps,\db$-boules pour $\eps>0$ petit.

Consid\'erons $B_P(u,\eps)$ avec $u$ extrait d'un bloc d'ordre quelconque avec $|u|\geq l_{n_0}$. Soit $n$ tel que
$l_n\leq |u|<l_{n+1}$.

\medbreak
$\underline{\text{\sl Premier cas:}}$ les symboles $0$ occupent au moins $99/100$e de $u$. Donc ils
occupent au moins $9/10$e de tout $v$ avec $\db(u,v)<9/100$ (et $|v|=|u|$).

Vu sa longueur, $v$ ne peut contenir de $n+1$-bloc complet. $v$ rencontre donc au plus un intervalle
intercalaire d'ordre $\geq n+1$. Le reste de $v$ est constitu\'e de z\'ero, un ou deux
intervalles. Chacun de ces intervalles est constitu\'e par la concat\'enation de $n-1$-blocs chacun
suivi de l'intervalle intercalaire d'ordre $n-1$ associ\'e ainsi que d'intervalles intercalaires
d'ordre $n$. En effet, on peut oublier les blocs ou intervalles incomplets, leur longueur
(major\'ee par $\max(l_{n-1},l_n^\beta)$) \'etant tr\`es petite devant $l_n$. La densit\'e des $0$ sur
ces intervalles est donc major\'ee par: 
 $$
    \mu_{n-1}(0) + \frac{1}{l_{n-1}^{1-\beta}} + \frac{1}{l_n^{1-\beta}} < 1/100.
 $$
En effet $\mu_{n-1}(0)$, la fraction des $B_{n-1}^i$ occup\'ee par $0$ est arbitrairement petite
et $l_{n-1},l_n$ sont arbitrairement grands, car $l_0$ est grand.

L'intervalle intercalaire d'ordre $n+1$ doit donc rencontrer $v$ sur une fraction $x$ de sa
longueur telle que $x+\frac{1}{100}(1-x)\geq\frac{9}{10}$. Donc $x\geq 89/99$. En particulier:
 $$
     v_{11/99|v|}\dots v_{89/99|v|}=\underbrace{0\dots0}_{79/99|v|}.
 $$
On voit donc que $B_P(u,9/100)\subset T^{-11/99|u|}([0^\ell]_P)$ avec $\ell\eqdef79/99|u|$.

$u$ correspond donc \`a un segment intercalaire dans un bloc d'ordre $\geq N$, avec $N$ minimal tel que
$l_N^\beta\geq\ell$:
 $$
    \mu([0^\ell]_P) \leq \sum_{n\geq N} \frac{l_n^\beta}{l_n}
        \leq \sum_{k\geq0} \frac{1}{(2^kl_N)^{1-\beta}}
        \leq \cte \cdot \frac{1}{|u|^{1/\beta-1}} \leq \frac{1}{|u|^{\gamma}}
 $$
car $l_{n+1}/l_n\geq2$ vu l'estimation sur $l_n$. La derni\`ere in\'egalit\'e a lieu d\`es que:
 $
    \beta<\frac{1}{\gamma+1}.
 $
On a donc bien: $\mu(B_P(u,\eps))\leq 1/|u|^\Gamma$ dans ce cas.

\medbreak
$\underline{\text{\sl Deuxi\`eme cas:}}$ la densit\'e des $0$ est inf\'erieure \`a $99/100$. En particulier, la
portion de $u$ occup\'ee par les un ou deux intervalles intercalaires d'ordre au moins $n+1$ rencontrant 
$u$ est major\'ee par $99/100$. Quitte \`a r\'eduire la longueur de $u$ par un facteur $100$, on peut donc supposer
que $u$ ne rencontre pas de tel intervalle. $u$ est donc inclus dans un $n+1$-bloc. Soit $\eps_*>0$ fourni 
par l'affirmation~(\clmSystNonPiste).

Supposons d'abord que $|u|\geq l'_{n+1}=l_{n+1}^\beta$. Alors $B_P(u,\eps_*)$ est inclus dans un \'etage bien
d\'efini de la tour associ\'ee au $n+1$-bloc contenant $u$. Donc:
 $$
    \mu(B_P(u,\eps_*))\leq \frac{1}{l_{n+1}N_{n+1}} = \frac{1}{l_{n+1}^{\gamma}}
        \leq \frac{1}{|u|^{\gamma}} \leq \frac{1}{|u|^\Gamma}.
 $$

Supposons maintenant que $l_n/100\leq |u|<l'_{n+1}=l_n^{1+\beta}$. On a:
 $$
    B_P(u,\eps_*) \subset \bigcup_{k=1}^{r}
       T^{-d_k} B_P(u_k,\eps_*)
 $$
o\`u $u=u_1\dots u_r$ est la d\'ecomposition de $u$ suivant les $n$-blocs et $d_k=|u_1\dots
u_{k-1}|$. On a $r\leq |u|/l_n+2\leq 2l_n^{\beta}$. Quitte \`a r\'eduire $u$ d'une fraction de sa
longueur major\'ee par $\frac{2l'_n}{l_n/200}<<1$, on suppose que chaque $u_k$ est de longueur au moins
$l'_n=l_n^\beta$. $B_P(u_k,\eps_*)$ est donc inclus dans un \'etage bien d\'efini de la tour correspondant \`a
un $n$-bloc bien d\'efini. Donc:
 $$
    \mu(B_P(u,\eps_*)) \leq 2l_n^{\beta} \cdot \frac{1}{l_nN_n}
       = \frac{2}{l_n^{\gamma-\beta}} \leq \frac{\cte}{|u|^{(\gamma-\beta)/(1+\beta)}}
       < \frac{1}{|u|^\Gamma}
 $$
d\`es que $\beta>0$ est assez petit pour que:
 $
    \frac{\gamma-\beta}{1+\beta} > \Gamma.
 $

\enddocument

\newpage

\bigbreak
\hskip-0.5cm\hfill \scaledpicture 165mm by 77mm (deftour scaled 744) \hfill\hbox{}
\centerline{\smc Figure 1.}
\centerline{\sl Exemple de tour irr\'eguli\`ere de hauteurs 3, 7 et 11.}
\bigbreak

\medbreak
\hskip-1cm\hfill \scaledpicture 167mm by 116mm (defemboite scaled 700) \hfill\hbox{}
\centerline{\smc Figure 2.}
\centerline{\sl Exemple d'embo\^{\i}tement de $\tau'$ dans $\tau$.}

\medbreak

\enddocument